\documentclass[sn-mathphys]{sn-jnl}

\jyear{2022}%

\theoremstyle{thmstyleone}%
\newtheorem{theorem}{Theorem}[section]
\theoremstyle{thmstyletwo}%
\newtheorem{remark}[theorem]{Remark}

\theoremstyle{thmstylethree}%

  \makeatother
  
  \makeatletter
\newcommand\fs@nocaptionruled{
  \let\@fs@capt\relax
  \def\@fs@pre{}
  \def\@fs@post{\kern2pt\hrule\relax}%
  \def\@fs@mid{\kern2pt\hrule\kern2pt}%
  \let\@fs@iftopcapt\iftrue}
\makeatother

\makeatletter
\algnewcommand{\LineComment}[1]{\Statex \hskip\ALG@thistlm \(\triangleright\) #1}
\makeatother
\usepackage{multirow,xfrac,array}
\usepackage{graphics}
\usepackage{subcaption}
\usepackage{algpseudocode}
\usepackage{eqparbox}
\usepackage{enumitem}

\newcommand{\alglinenoNew}[1]{\newcounter{ALG@line@#1}}

\newcommand{\alglinenoPush}[1]{\setcounter{ALG@line@#1}{\value{ALG@line}}}

\newcommand{\bu}{\pmb{u}}
\newcommand{\bv}{\pmb{v}}
\newcommand{\bn}{\pmb{n}}
\newcommand{\bK}{\pmb{K}}
\newcommand{\bP}{\pmb{P}}
\newcommand{\iS}{\mathcal{S}}
\newcommand{\iP}{\mathcal{P}}
\newcommand{\iD}{\mathcal{D}}
\newcommand{\iF}{\mathcal{F}}
\newcommand{\iT}{\mathcal{T}}
\newcommand{\Rv}[1]{\textcolor{black}{#1}}
\begin{document}

\title[Domain Decomposition for Reduced Fracture Models]{Fast and accurate domain decomposition methods for reduced fracture models with nonconforming time grids}

\author[1]{\fnm{} \sur{Phuoc-Toan Huynh}}\email{tph0017@auburn.edu}
\author[1]{\fnm{} \sur{Yanzhao Cao}}\email{yzc0009@auburn.edu}
\author*[1]{\fnm{} \sur{Thi-Thao-Phuong~Hoang}}\email{tzh0059@auburn.edu}

\affil[1]{\orgdiv{Department of Mathematics and Statistics}, \orgname{Auburn University}\orgaddress{\street{}, \city{Auburn}, \postcode{36849}, \state{Alabama}, \country{USA}}}

\abstract{This paper is concerned with the numerical solution of compressible fluid flow in a fractured porous medium.  The fracture represents a fast pathway (i.e., with high permeability) and is modeled as a hypersurface embedded in the porous medium.  We aim to develop fast-convergent and accurate global-in-time domain decomposition (DD) methods for such a reduced fracture model,  in which smaller time step sizes in the fracture can be coupled with larger time step sizes in the subdomains.  Using the pressure continuity equation and the tangential PDEs in the fracture-interface as transmission conditions,  three different DD formulations are derived; each method leads to a space-time interface problem which is solved iteratively and globally in time.  Efficient preconditioners are designed to accelerate the convergence of the iterative methods while preserving the accuracy in time with nonconforming grids.  Numerical results for two-dimensional problems with non-immersed and partially immersed fractures are presented to show the improved performance of the proposed methods. }




\maketitle

\section{Introduction}\label{sec1}
Numerical simulation of flow and transport in a fractured porous medium is challenging due to the presence of multiple spatial and temporal scales and the strong physical property heterogeneity of the domain of calculation.  In particular,  a fracture can represent either a fast pathway or a geological barrier,  depending on whether its permeability is much higher or much lower than the surrounding rock matrix.  Thus the time scales in the fractures and in the rock matrix may vary significantly.  In addition,  the width of the fracture is much smaller than the size of the domain of calculation and any reasonable spatial mesh size.  To avoid local refinement around the fractures, one possible approach is to reduce the original problem into a new one where the fractures are treated as domains of co-dimension one, i.e, interfaces between subdomains (see \cite{1, 2, 32,4, 12, 40,Gander2021,SHLee} and the references therein).  Models with such low-dimensional fractures are known as reduced fracture models or mixed-dimensional models. 

In this paper, we are concerned with numerical algorithms for a reduced fracture model of compressible fluid flow in which the fracture has larger permeability than the surrounding porous medium.  For such a case,  the fluid flows rapidly through the fracture while it moves much more slowly through the rock matrix.  Hence,  using a single-time step size throughout the entire domain of calculation is computationally inefficient.  This work aims to develop fast-convergent and accurate global-in-time domain decomposition (DD) methods for the reduced fracture model in which smaller time step sizes in the fracture can be coupled with larger time step sizes in the subdomains. For the spatial discretization of the flow problem, we use mixed finite elements as they are mass conservative and can handle well heterogeneous and anisotropic diffusion tensors \cite{7, 44}.  

Global-in-time DD methods provide a powerful tool to perform parallel simulations of time-dependent physical phenomena with different time steps across the domain.  
These methods are obtained by decoupling the given dynamic system into dynamic subsystems defined on the subdomains (resulting from a spatial decomposition), then time-dependent problems are solved in each subdomain at each iteration and information is exchanged over space-time interfaces between subdomains.  
Global-in-time DD is different from the classical DD approach~\cite{41, 42} where the model problem is first discretized in time by an implicit scheme, then at each time step the iteration is performed and involves the solution of stationary problems in the subdomains.  The same time step is required for the classical approach, while for global-in-time DD,  local time discretizations can be enforced in different regions of the domain. 

There are basically two types of global-in-time DD methods. The first type is based on the physical transmission conditions, for example, the Dirichlet-Neumann and Neumann-Neumann waveform relaxation methods~\cite{34, 33, 18, 19,30}. The second type is based on more general transmission conditions such as Robin or Ventcel~\cite{Ventcel} conditions.  An important class of methods in this category is called the Optimized Schwarz Waveform Relaxation (OSWR) algorithm~\cite{14, 15, 5, 21,  23, 6} where additional coefficients involved in the transmission conditions are optimized to improve convergence rates.  Both approaches were used with mixed formulations to treat the pure diffusion problem in \cite{24} and the linear advection-diffusion problem in~\cite{27}. In particular, the global-in-time primal Schur (GTP-Schur) and global-in-time optimized Schwarz (GTO-Schwarz) methods were proposed in \cite{24,27}.  For each method,  an interface problem on the space-time interfaces between subdomains is derived and is solved iteratively over the whole time interval.  

In \cite{26},  GTP-Schur and GTO-Schwarz methods were studied for a reduced fracture model of a single-phase, compressible fluid flow in a porous medium with a ``fast-path'' fracture.  For such a model,  the physical transmission conditions consist of the pressure continuity equation and the tangential PDEs in the fracture.  Based on these conditions,  a space-time interface problem for GTP-Schur is obtained using the time-dependent Dirichlet-to-Neumann operator.  Two preconditioners were considered in \cite{26}: the local preconditioner and the time-dependent Neumann-Neumann preconditioner.  The former is adapted from \cite{3} (for second-order elliptic PDEs) and the latter is an extension of the balancing domain decomposition (BDD) preconditioner~\cite{35, 9, 36} to time-dependent problems. 
The GTO-Schwarz method uses the so-called Ventcel-to-Robin transmission conditions which are obtained by taking linear combinations of the pressure continuity equation and the PDEs in the fracture.  
These new transmission conditions contain a
free parameter, which is used to accelerate the convergence of the iterative method.   The interface problem for GTO-Schwarz is derived using the Ventcel-to-Robin operator,  and requires no preconditioner.  Different time steps in the fracture and in the rock matrix can be used for both GTP-Schur and GTO-Schwarz via a suitable $L^2$ projection in time.  An optimal projection algorithm can be found in \cite{16, 17}. 
 
The global-in-time DD methods proposed in \cite{26} have two drawbacks.  Firstly, the preconditioners for GTP-Schur are not effective: numerical results in~\cite{26} show that the convergence of GTP-Schur with either local or Neumann-Neumann preconditioner is much slower than that of GTO-Schwarz.  Secondly,  while GTO-Schwarz converges remarkably fast,  it does not preserve the accuracy in time in the fracture with nonconforming time grids.   In particular,  using a smaller time step in the fracture than in the surrounding rock matrix does not improve the errors in the fracture,  compared to using the same time step in the whole domain.  This is also the case for GTP-Schur with the Neumann-Neumann preconditioner. 
%

In this paper, we develop efficient global-in-time DD methods, based on physical transmission conditions,  which overcome the difficulties encountered in~\cite{26}.  
The contributions of this work include four aspects. Firstly,  an efficient preconditioner is derived to enhance the convergence of GTP-Schur.  The new preconditioner, namely Ventcel-Ventcel preconditioner,  provides a more accurate approximation of the (pseudo) inverse of the interface operator associated with the GTP-Schur method.  The Ventcel-Ventcel preconditioner involves solving the subdomain problems with Ventcel boundary conditions, instead of Neumann conditions as for the Neumann-Neumann preconditioner in \cite{26}. Secondly,  we introduce the global-in-time dual Schur (GTD-Schur) method in which the interface problem is derived using the time-dependent Neumann-to-Dirichlet operator, instead of the Dirichlet-to-Neumann operator as for the primal Schur approach. 
The dual formulation was first proposed for mixed finite elements~\cite{20} and late on widely studied for finite elements in finite element tearing and interconnecting (FETI) methods \cite{Farhat1991, Farhat1995}. To the best of our knowledge,  the global-in-time dual Schur approach has not been use yet to study the reduced fracture models in the literature.  We also introduce the so-called Dirichlet-Dirichlet preconditioner for this method to enhance its performance.  Thirdly,  we propose a new method, namely global-in-time fracture-based Schur (GTF-Schur), by combining the ideas of the primal and dual Schur methods.  One advantage of this new method is that the space-time interface operator is close to the identity operator; as a consequence, the iterative solver works well without requiring any preconditioners.  Lastly, we carry out numerical experiments for both non-immersed and partially immersed fractures to verify and compare the performance of the proposed methods with different time steps in the fracture and in the rock matrix. 
 
The rest of this paper is organized as follows: in Section~\ref{sec2} we present the model problem in mixed form and describe briefly the reduced process to transform the original problem into the reduced fracture model.  The Ventcel-Ventcel preconditioner for the GTP-Schur method is constructed in Section~\ref{sec3}.  In Section~\ref{sec4}, we formulate the GTD-Schur method and its Dirichlet-Dirichlet preconditioner. The GTF-Schur method is developed in Section~\ref{sec5}. In Section~\ref{sec6}, the semidiscrete problems for all proposed methods in time using different time grids in the subdomains are considered.  Numerical results are presented in Section~\ref{sec7} to illustrate and compare the performance of the proposed methods with GTO-Schwarz.  The paper is then closed with a conclusion section.
\section{A reduced fracture model}\label{sec2}
Let $\Omega$ be a bounded domain in $\mathbb{R}^{d}$ $(d= 2, 3)$ with Lipschitz boundary $\partial\Omega$, and $T>0$ be some fixed time.  Consider the flow problem of a single phase, compressible fluid written in mixed form as follows:
\begin{equation}
\label{original_problem}
\begin{array}{clll}
\phi\partial_t{p} + \text{div } \bu & = & q &\text{ in } \Omega \times (0, T), \\
\bu & = & -\bK \nabla{p} & \text{ in } \Omega \times (0, T), \\
p & = & 0 & \text{ on } \partial\Omega \times (0, T), \\
p(\cdot, 0) & = & p_0 & \text{ in } \Omega,
\end{array}
\end{equation}
where $p$ is the pressure, $\bu$ the velocity, $q$ the source term, $\phi$ the storage coefficient, and $\bK$ a symmetric, time-independent, hydraulic, conductivity tensor.
Suppose that the fracture $\Omega_f$ is a subdomain of $\Omega$, whose thickness is $\delta$, that separates $\Omega$ into two connected subdomains: $\Omega \backslash \overline{\Omega}_f = \Omega_1 \cup \Omega_2, $ and $\Omega_1 \cap \Omega_2 = \emptyset.$
For simplicity, we assume further that $\Omega_f$ can be expressed as \vspace{-0.2cm}
$$
\Omega_f  = \left\{ \textbf{\textit{x}} \in \Omega: \textbf{\textit{x}} = \textbf{\textit{x}}_{\gamma} + s\bn \text{ where } \textbf{\textit{x}}_{\gamma} \in \gamma \text{ and } s \in \left(-\dfrac{\delta}{2}, \dfrac{\delta}{2}\right)\right\}, \vspace{-0.2cm}
$$
where $\gamma$ is the intersection between a line $(d=2)$ or a plane $(d=3)$ with $\Omega$.
\begin{figure}[H]
\vspace{-0.4cm}
\centering
\includegraphics[scale=0.55]{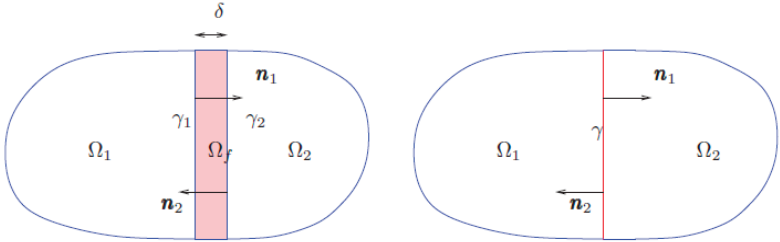}
\caption{The domain $\Omega$ with the fracture $\Omega_{f}$ (left) and the fracture-interface~$\gamma$~(right).}
\label{nonconforming_time}  \vspace{-0.3cm}
\end{figure}
We denote by $\gamma_i$ the part of the boundary of $\Omega_i$ shared with the boundary of the fracture $\Omega_f$: $\, \gamma_i = (\partial\Omega_i \cap \partial\Omega_f) \cap \Omega$,  for $i =1,2$.  Let $\bn_i$ be the unit, outward pointing, normal vector field on $\partial\Omega_i$, where $\bn= \bn_1 = - \bn_2$. For $i=1,\; 2,\; f$, and for any scalar, vector, or tensor valued function $\phi$ defined on $\Omega$, we denote by $\phi_i$ the restriction of $\phi$ to $\Omega_i$. The original problem \eqref{original_problem} can be rewritten as the following transmission problem: \vspace{-0.2cm}
\begin{equation}
\label{multidomain_problem}
\begin{array}{cllll}
\phi_i\partial_t{p_i} + \text{div }\bu_i & = & q_{i} & \text{ in } \Omega_i \times (0, T), & i = 1, 2, f, \\
\bu_i & = & - \bK_i \nabla{p}_i & \text{ in } \Omega_i \times (0, T), & i = 1, 2, f, \\
p_i & = & 0 & \text{ on } \left(\partial\Omega_i \cap \partial\Omega\right) \times (0, T), & i=1, 2, f, \\
p_i & = & p_f & \text{ on } \gamma_i \times (0, T), & i=1, 2, \\
\bu_i\cdot \bn_i & = & \bu_f\cdot\bn_i & \text{ on } \gamma_i \times (0, T), & i=1, 2, \\
p_i(\cdot, 0) & = & p_{0, i} & \text{ in } \Omega_i, & i=1, 2, f.
\end{array} \vspace{-0.2cm}
\end{equation}
The reduced fracture model that we consider in this paper was first proposed in~\cite{1, 2} under the assumption that the fracture has larger permeability than that in the rock matrix. The model is obtained by averaging across the transversal cross sections of the $d$-dimensional fracture $\Omega_f$.  
We use the notation $\nabla_{\tau}$ and $\text{div}_{\tau}$
for the tangential gradient and tangential divergence, respectively. \Rv{We write $\phi_{\gamma}$ and $\bK_{\gamma}$ for $\delta\phi_f$ and $\bK_{f, \tau}$, respectively, where $\bK_{f, \tau}$ is the tangential component of $\bK_f$}. The reduced model consists of equations in the subdomains, \vspace{-0.2cm}
\begin{equation}
\label{reduced_subdomain}
\left.\begin{array}{rcll}
\phi_i\partial_t{p_i}+\text{div }\bu_i&=&q_{i} &\text{ in } \Omega_i\times (0, T), \\
\bu_i&=&-\bK_i\nabla{p_i} &\text{ in } \Omega_i\times (0, T), \\
p_i&=&0 &\text{ on } \left(\partial\Omega_i \cap \partial\Omega\right) \times (0, T), \\
p_i&=&p_{\gamma} &\text{ on } \gamma \times (0, T), \\
p_i(\cdot, 0)&=&p_{0, i} &\text{ in } \Omega_i, 
\end{array}\right.  \vspace{-0.2cm}
\end{equation}
for $i=1,2,$ and equations in the fracture, \vspace{-0.2cm}
\begin{equation}
\label{reduced_fracture}
\begin{array}{rcll}
\phi_{\gamma}\partial_t{p_{\gamma}}+\text{div}_{\tau }\bu_{\gamma}&=& q_{\gamma} +\sum^{2}_{i=1}\left( \bu_i \cdot \bn_i\right)_{\vert \gamma} & \text{ in } \gamma \times (0, T), \\
\bu_{\gamma} &=&-\bK_{\gamma}\delta\nabla_{\tau}p_{\gamma} & \text{ in } \gamma \times (0, T), \\
p_{\gamma}&=&0 &\text{ on } \partial\gamma \times (0, T), \\
p_{\gamma}(\cdot, 0)&=&p_{0, \gamma} & \text{ in } \gamma.
\end{array} \vspace{-0.2cm}
\end{equation}
To write the weak formulation of \eqref{reduced_subdomain}-\eqref{reduced_fracture}, we use the convention that if $V$ is a space of functions, then $\pmb{V}$ is a space of vector functions having each component in $V$.  For arbitrary domain $\mathcal{O}$, we denote by $\left(\cdot, \cdot\right)_{\mathcal{O}}$ the inner product in $L^2\left(\mathcal{O}\right)$ or $\mathbf{\textbf{\textit{L}}^2\left(\mathcal{O}\right)}$. 
We next define the following Hilbert spaces: \vspace{-0.2cm}
\begin{align*}
\begin{array}{rl}
M &= \left\{v  = \left( v_1, v_2, v_{\gamma}\right) \in L^2\left(\Omega_1\right) \times L^2\left(\Omega_2\right) \times L^2\left(\gamma\right) \right\},  \vspace{0.2cm}\\
\Sigma &= \left\{\bv = \left(\bv_1, \bv_2, \bv_{\gamma}\right) \in \textbf{\textit{L}}^2\left(\Omega_1\right) \times \textbf{\textit{L}}^2\left(\Omega_2\right) \times \textbf{\textit{L}}^2\left(\gamma\right): \; \text{div} \; \bv_i \in L^2\left(\Omega_i\right), \; i=1,2, \right. \\
& \qquad \left.\text{and } \text{div}_{\tau} \; \bv_{\gamma} - \sum\limits^2_{i=1} \bv_i\cdot\bn_{i \vert \gamma} \in L^2(\gamma)\right\}.
\end{array}
\end{align*}
We define the bilinear forms $a(\cdot, \cdot)$, $b(\cdot, \cdot)$ and $c(\cdot, \cdot)$ on $\Sigma \times \Sigma$,  $\Sigma \times M$, and $M \times M$, respectively, and the linear form $L_q$ on $M$ by \vspace{-0.2cm}
\begin{align*}
a\left(\bu, \bv\right) &= \sum\limits^2_{i=1} \left(\bK^{-1}_i\bu_i, \bv_i\right)_{\Omega_i} + \left(\left(\bK_{\gamma}\delta\right)^{-1}\bu_{\gamma}, \bv_{\gamma}\right)_{\gamma}, \\
 b\left(\bu, \mu\right) &= \sum\limits^2_{i=1}\left(\text{div}\; \bu_i, \mu_i\right)_{\Omega_i} + \left(\text{div}_{\tau} \; \bu_{\gamma} - \sum\limits^2_{i=1} \bu_i\cdot\bn_{i \vert \gamma}, \mu_{\gamma}\right)_{\gamma},\\
c_{\phi}(\eta, \mu) &= \sum\limits^2_{i=1}\left(\phi_i\eta_i, \mu_i\right)_{\Omega_i} + \left(\phi_{\gamma}\eta_{\gamma}, \mu_{\gamma}\right)_{\gamma}, \quad L_{q}(\mu) = \sum\limits^2_{i=1}(q_i, \mu_i)_{\Omega_i}. \vspace{-0.2cm}
\end{align*}
%
%
The weak form of \eqref{reduced_subdomain}-\eqref{reduced_fracture} can be written as follows: 

Find $p \in H^1(0, T; M)$ and $\bu \in L^2(0, T;\Sigma)$ such that \vspace{-0.2cm}
\begin{equation}
\label{weak_reduced}
\begin{array}{rcll}
a\left(\bu, \bv\right) - b\left(\bv, p \right) & = & 0 & \forall \bv \in \Sigma, \\
c_{\phi}\left(\partial_{t}p, \mu\right) + b\left(\bu, \mu\right) &=& L_q(\mu) & \forall \mu \in M,
\end{array} \vspace{-0.2cm}
\end{equation}
together with the initial conditions: \vspace{-0.2cm}
\begin{equation}
\label{initial_weak_reduced}
p_i(\cdot, 0)  =  p_{0, i}, \; \text{in} \; \Omega_i, \; i=1,2,  \quad \text{and} \quad p_{\gamma}(\cdot, 0)  =  p_{0, \gamma}, \; \text{in} \; \gamma.\vspace{-0.2cm}
\end{equation}
The well-posedness of problem \eqref{weak_reduced}-\eqref{initial_weak_reduced} is given by the following theorem. The reader is referred to \cite[Theorem $2.1$]{26} for the details of the proof.  \vspace{-0.3cm}
\begin{theorem}{\cite{26}} 
Assume that the storage coefficient $\phi_{i}$, $i=1,2,\gamma$ is bounded above and below by positive constants, and that there exist positive constants $\bK_{-}$ and $\bK_{+}$ such that
\begin{itemize}
\item[(i)]$\zeta^T\bK^{-1}_i(x) \zeta \geq \bK_{-}\vert \zeta \vert^2,$ and $\vert \bK_i(x)\zeta \vert \leq \bK_{+} \vert \zeta \vert,$ for a.e.  $x \in \Omega_i$ and $ \forall \zeta \in \mathbb{R}^d, \; i=1,2,$
\item[(ii)] $\eta^T\left(\bK_{\gamma}(x)\delta\right)^{-1}\eta \geq \bK_{-}\vert \eta \vert^2$ and $\vert \left(\bK_{\gamma}(x)\delta\right)^{-1}\eta\vert \leq \bK_{+}\vert \eta \vert$, for a.e.  $x \in \gamma$ and $\forall \zeta \in \mathbb{R}^{d-1}$.
\end{itemize}
Given $q$ in $L^2(0, T; M)$ and $p_0 =\left(p_{0, 1}, \; p_{0, 2}, \; p_{0, \gamma}\right)$ in $H^1_{*}$,  where \vspace{-0.1cm}
\begin{align*}
H^1_{*} := \left\{\mu = \left(\mu_1, \mu_2, \mu_{\gamma}\right) \in H^1\left(\Omega_1\right) \times H^1\left(\Omega_2\right) \times H^1_0\left(\gamma\right): \; \mu_i = 0 \text{ on } \partial\Omega_i \cap \partial\Omega, \right.\vspace{-0.2cm}\\
\qquad \qquad \quad \left. \text{ and } \mu_i - \mu_{\gamma} =0 \text{ on } \gamma, \; i=1, 2\right\}. \vspace{-0.2cm}
\end{align*}
Then problem  \eqref{weak_reduced}-\eqref{initial_weak_reduced} has a unique solution $\left(p, \bu\right) \in H^1(0, T;M) \times L^2\left(0, T; \Sigma\right).$ \vspace{-0.4cm}
\end{theorem}
%
We shall use global-in-time DD to find a numerical solution of problem \eqref{weak_reduced}-\eqref{initial_weak_reduced} with different time steps in the fracture and the surrounding medium.  The DD formulation can be obtained by treating the fractures as a (physical) interface between subdomains with the following transmission conditions: \vspace{-0.2cm}
\begin{align}
&\hspace{2.3cm} p_i=p_{\gamma}, \hspace{3cm} \text{ on } \gamma \times (0, T), \label{transmission_1} \\
&\begin{array}{rcll}
\phi_{\gamma}\partial_t{p_{\gamma}}+\text{div}_{\tau }\bu_{\gamma}&=& q_{\gamma} +\sum^{2}_{i=1}\left( \bu_i \cdot \bn_i\right)_{\vert \gamma} & \text{ in } \gamma \times (0, T), \\
\bu_{\gamma} &=&-\bK_{\gamma}\delta\nabla_{\tau}p_{\gamma} & \text{ in } \gamma \times (0, T), \\
p_{\gamma}&=&0 &\text{ on } \partial\gamma \times (0, T), \\
p_{\gamma}(\cdot, 0)&=&p_{0, \gamma} & \text{ in } \gamma.
\end{array} 
\label{transmission_2} \vspace{-0.2cm}
\end{align}
In the next sections, three global-in-time DD methods are derived based on these physical transmission conditions.  For each method, a space-time interface problem is formulated and solved iteratively.  
%
\section{Global-in-time primal Schur (GTP-Schur) method}\label{sec3}
The idea of GTP-Schur is to impose \eqref{transmission_1} as Dirichlet boundary conditions for the subdomain problems: \vspace{-0.2cm}
\begin{equation}
p_{i}=\lambda, \quad \text{on} \; \gamma \times (0,T), \; i=1,2, \vspace{-0.2cm}
\end{equation}
where $\lambda$ represents the fracture pressure $p_{\gamma}$. Then a space-time interface problem is formed by enforcing the remaining transmission condition \eqref{transmission_2}.  To derive the formulation of GTP-Schur,  we define the Dirichlet-to-Neumann operators $\iS^{\text{DtN}}_i, \; i=1, 2$: \vspace{-0.2cm}
\begin{align*} 
\small
\begin{array}{ccc}
\iS^{\text{DtN}}_{i}: H^1(0, T;H^{\frac{1}{2}}_{00}(\gamma)) \times L^2(0, T; L^2\left(\Omega_i\right)) \times H^{1}_{*, \gamma}\left(\Omega_i\right) & \longrightarrow & L^2 \left (0, T; (H^{\frac{1}{2}}_{00}(\gamma))^{\prime}\right ), \\
\left(\lambda, q_{i}, p_{0, i}\right) & \longmapsto & \bu_i\cdot\bn_{i\vert\gamma},
\end{array} \label{eq:S_DtN}\vspace{-0.2cm}
\end{align*}
where $H^{1}_{*, \gamma}\left(\Omega_i\right) := \left\{ \mu \in H^{1}\left(\Omega_i\right): \mu = 0 \text{ on } \left(\partial\Omega_i \cap \partial\Omega\right)\right\}$ and $(p_i, \bu_i)$ is the solution of the problem \vspace{-0.2cm}
\begin{equation}
\label{subdom_solver_Schur}
\left.\begin{array}{rcll}
\phi_i\partial_t{p_i}+\text{div }\bu_i&=&q_{i} &\text{ in } \Omega_i\times (0, T), \vspace{2pt}\\
\bu_i&=&-\bK_i\nabla{p_i} &\text{ in } \Omega_i\times (0, T), \\
p_i&=&0 &\text{ on } \left(\partial\Omega_i \cap \partial\Omega\right) \times (0, T), \\
p_i&=&\lambda &\text{ on } \gamma \times (0, T), \\
p_i(\cdot, 0)&=&p_{0, i} &\text{ in } \Omega_i.
\end{array}\right. \vspace{-0.2cm}
\end{equation}
The space-time interface problem with unknown $\lambda$ reads as: \vspace{-0.2cm}
\begin{equation}
\label{interface_problem_Schur_a}
\begin{array}{rcll}
\phi_{\gamma}\partial_t\lambda+\text{div}_{\tau }\bu_{\gamma}&=& q_{\gamma} +\sum^{2}_{i=1}\iS^{\text{DtN}}_{i}(\lambda, q_{i}, p_{0, i}) & \text{ in } \gamma \times (0, T), \\
\bu_{\gamma} &=&-\bK_{\gamma}\delta\nabla_{\tau}\lambda & \text{ in } \gamma \times (0, T), \\
\lambda&=&0 &\text{ on } \partial\gamma \times (0, T), \\
\lambda(\cdot, 0)&=&p_{0, \gamma} & \text{ in } \gamma,
\end{array} \vspace{-0.2cm}
\end{equation}
or equivalently, 
\begin{equation}
\label{interface_problem_Schur_b}
\small
\begin{array}{rcll}
\phi_{\gamma}\partial_t\lambda+\text{div}_{\tau }\bu_{\gamma} -\sum^{2}_{i=1}\iS^{\text{DtN}}_{i}(\lambda, 0,0)&=& q_{\gamma} & \hspace{-1.2cm}+\sum^{2}_{i=1}\iS^{\text{DtN}}_{i}(0, q_{i}, p_{0, i}) \vspace{0.1cm}\\
&& &\text{ in } \gamma \times (0, T), \\
\bu_{\gamma} &=&-\bK_{\gamma}\delta\nabla_{\tau}\lambda & \text{ in } \gamma \times (0, T), \\
\lambda&=&0 &\text{ on } \partial\gamma \times (0, T), \\
\lambda(\cdot, 0)&=&p_{0, \gamma} & \text{ in } \gamma,
\end{array} \vspace{-0.1cm}
\end{equation} 
or in compact form (space-time), \vspace{-0.2cm}
\begin{equation}
\label{interface_compact_Schur}
\iS_{\iP}\begin{pmatrix}
\lambda 
\end{pmatrix} = \chi_{\iP}. \vspace{-0.2cm}
\end{equation}
Note that from the second equation of \eqref{interface_problem_Schur_b}, $\bu_{\gamma}$ is a function in $\lambda$, hence, the right-hand side operator of \eqref{interface_compact_Schur} is actually an operator in  only one variable $\lambda$.  

The space-time problem \eqref{interface_compact_Schur} is solved iteratively using, e.g., GMRES.  \Rv{The resulting algorithm is matrix free as the discrete counterpart of $\iS_{\iP}$ is not computed explicitly.  At each GMRES iteration,  $\iS_{\iP}(\lambda)$ is obtained by first solving the subdomain problems~\eqref{subdom_solver_Schur} over the whole time interval, then using the tangential PDEs~\eqref{interface_problem_Schur_b} in the fracture-interface.}
The convergence of the iterative algorithm is known to be significantly slow,  thus finding a suitable preconditioner is necessary to accelerate the iteration.  Two preconditioners were introduced in~\cite{26}.  The local preconditioner, ${\bP}^{-1}_{\text{loc}} $,  is computed by finding the discrete counterpart of the operator $\left(\text{div}_{\tau}\left(\bK_{\gamma}\delta\nabla_{\tau}\right)\right)^{-1}$. This preconditioner was proposed first in \cite{3} for stationary problems using the fact that the second order operator $\left(\text{div}_{\tau} \; \left(\bK_{\gamma}\delta\nabla_{\tau}\right)\right)$ is the dominant term in the interface problem. 
The second preconditioner is the (time-dependent) Neumann-Neumann preconditioner, ${\bP}^{-1}_{\text{NN}} $,  obtained by computing the (pseudo-)inverse of the Dirichlet-to-Neumann operators.  Such a preconditioner involves the solution of the subdomain problems with Neumann boundary conditions on the fracture-interface.  For the case with no fracture, the Neumann-Neumann preconditioner has been shown to be effective~\cite{24}. 
However, for the considered reduced fracture model,  it has been shown numerically in \cite{26} that the convergence speed of the iterative algorithm combined with these preconditioners is still slow and not efficient, especially the local preconditioner. From the derivation of these preconditioners, it can be seen that they do not provide good approximations of the inverse of the space-time operator on the left-hand side of the first equation in \eqref{interface_problem_Schur_b}.  Based on this observation,  we derive a new preconditioner, namely the Ventcel-Ventcel preconditioner, in the following. 
%
\subsection*{Ventcel-Ventcel preconditioner}\label{subsec3}
As $\lambda$ represents the fracture pressure $p_{\gamma}$ and by the definition of $\iS^{\text{DtN}}_{i}$,  the left-hand side of the first equation in \eqref{interface_problem_Schur_b} can be rewritten as \vspace{-0.2cm}
\begin{equation}
\label{Ventcell_BC}
\phi_{\gamma}\partial_t{\lambda} + \text{div}_{\tau} \; \bu_{\gamma} - \sum\limits^2_{i=1} \iS^{\text{DtN}}_i\left(p_{\gamma}, 0, 0\right) = \phi_{\gamma}\partial_t{p_{\gamma}} + \text{div}_{\tau} \; \bu_{\gamma} - \sum\limits^2_{i=1}\bu_i\cdot\bn_{i\vert\gamma}. \vspace{-0.2cm}
\end{equation}
The right-hand side of this equation resembles Ventcel boundary conditions~\cite{26}.  Thus, the preconditioned system for \eqref{interface_problem_Schur_b} should be computed by solving the subdomain problems with such Ventcel boundary conditions \eqref{Ventcell_BC} (instead of with Neumann conditions as used for the Neumann-Neumann preconditioner).  To formulate local problems with Ventcel conditions,  we introduce the Lagrange multipliers $p_{i, \gamma}, \; i=1,2,$ with $p_{i, \gamma}$ representing the trace on the interface $\gamma$ of the pressure $p_i$ in the subdomain $\Omega_i$.  It follows from the continuity of the pressure across the interface that  \vspace{-0.2cm}
\begin{equation} \label{eq:pi_g}
p_{1, \gamma} = p_{2, \gamma} = p_{\gamma}, \text{ in } \gamma \times (0, T). \vspace{-0.2cm}
\end{equation}
We write the Darcy equation associated with each $p_{i, \gamma}$ in the fracture as \vspace{-0.2cm}
\begin{equation} \label{eq:ui_g}
\bu_{\gamma, i} := -\bK_\gamma \delta\nabla_{\tau}p_{i, \gamma}, \; \; \text{in} \; \gamma \times (0, T), \; i=1, 2. \vspace{-0.2cm}
\end{equation}
Note that $\bu_{\gamma, i}, \; i= 1, 2$ represents the tangential velocity in the fracture associated with the pressure $p_{i, \gamma}$, and $\bu_{\gamma, 1} = \bu_{\gamma, 2} = \bu_{\gamma}, \;\;  \text{in}\; \gamma \times (0, T)$ according to \eqref{eq:pi_g} and \eqref{eq:ui_g}.  With such notation, the subdomain problem with Ventcel boundary condition reads as: \vspace{-0.2cm}
\begin{equation}
\label{subdom_solver_Schur_preconditioner}
\hspace{-0.2cm}\begin{array}{rcll}
\phi_i\partial_t{p_i}+\text{div }\bu_i&=&0 &\text{ in } \Omega_i\times (0, T), \\
\bu_i&=&-\bK_i\nabla{p_i} &\text{ in } \Omega_i\times (0, T), \\
p_i&=&0 &\text{ on } \left(\partial\Omega_i \cap \partial\Omega\right) \times (0, T), \\
\phi_{\gamma}\partial_t{p_{i, \gamma}} + \text{div}_{\tau}\bu_{\gamma, i} -\bu_i\cdot\bn_{i\vert\gamma}&=&\theta &\text{ on } \gamma \times (0, T), \\
\bu_{\gamma, i} &= & -\bK_i\delta\nabla_{\tau}p_{i, \gamma} &\text{ on } \gamma \times (0, T), \\
p_{i, \gamma} & = & 0 & \text{ on } \partial\gamma \times (0, T), \\
p_i(\cdot, 0)&=&0 &\text{ in } \Omega_i,
\end{array}\vspace{-0.2cm}
\end{equation}
for $i=1, 2$, where $\theta$ is given Ventcel data.  It can be shown that problem~\eqref{subdom_solver_Schur_preconditioner} has a unique weak solution; interested readers are referred to \cite[Theorem 4.1]{26} for more details of the proof. 
Next, we define the following Ventcel-to-Dirichlet operator ${\iS}^{\text{VtD}}_{i}, \; i=1, 2$: \vspace{-0.2cm}
\begin{align*} 
\begin{array}{ccc}
\iS^{\text{VtD}}_{i}: L^2\left(0, T; L^2\left(\gamma\right)\right) & \longrightarrow & H^{1}\left(0, T; L^2\left(\gamma\right)\right), \\
\theta & \longrightarrow & p_{i, \gamma}, \vspace{-0.2cm}
\end{array}
\end{align*}
where $\left(p_i, \; \bu_i, \; p_{i, \gamma}, \; \bu_{\gamma, i}\right)$, $i=1, 2,$ is the solution of the subdomain problem~\eqref{subdom_solver_Schur_preconditioner}.  Then the Ventcel-Ventcel preconditioner $P^{-1}_{\text{VV}}$ for problem \eqref{interface_compact_Schur} is given by \vspace{-0.2cm}
\begin{align*}
\bP^{-1}_{\text{VV}} := \sigma_1{\iS}^{\text{VtD}}_{1} + \sigma_2{\iS}^{\text{VtD}}_{2}, \vspace{-0.2cm}
\end{align*}
where $\sigma_i : \gamma \times (0, T) \rightarrow [0, 1]$ is such that $\sigma_1 + \sigma_2=1$. The preconditioned system for \eqref{interface_compact_Schur} with the Ventcel-Ventcel preconditioner is defined as: \vspace{-0.2cm}
\begin{equation}
\label{preconditioner_system_Ventcell}
\bP^{-1}_{\text{VV}}\left(\iS_{\iP}(\lambda)\right) = P^{-1}_{\text{VV}}(\chi_{\iP}), \; \; \text{in } \gamma \times (0, T). \vspace{-0.2cm}
\end{equation}
\indent \Rv{We summarize the GTP-Schur method with the Ventcel-Ventcel preconditioner in Algorithm~\ref{GTP_Schur}.  Note that the operator $\bP^{-1}_{\text{VV}}$ can be replaced by $\bP^{-1}_{\text{NN}}$ (i.e., the Neumann-Neumann preconditioner) or by the identity operator (i.e., no preconditioner).  We will compare numerical performance of these algorithms and verify the improvement by the Ventcel-Ventcel preconditioner in Section~\ref{sec7}. }  \vspace{-0.2cm} 
\begin{algorithm}[H]
\caption{GTP-Schur method with Ventcel-Ventcel preconditioner}
\label{GTP_Schur}
\textbf{Input:} initial guess $\lambda^{(0)}$, stopping tolerance $ 0 <\epsilon \ll 1$, maximum number of iterations $N_{\text{max}}$. \\
\textbf{Output:} space-time fracture pressure $\lambda$. 
\begin{algorithmic}[1]
\State Compute $\chi_{\iP} = q_{\gamma} +\sum^{2}_{i=1}\iS^{\text{DtN}}_{i}(0, q_{i}, p_{0, i}).$
\State Evaluate $\iS_{\iP}(\lambda^{(0)}) = \phi_{\gamma}\partial_t\lambda^{(0)} -\text{div}_{\tau }\bK_{\gamma}\delta\nabla_{\tau}\lambda^{(0)} -\sum\limits^{2}_{i=1}\iS^{\text{DtN}}_{i}(\lambda^{(0)}, 0,0).$
\State Set $r_0 := \chi_{\iP} - \iS_{\iP}(\lambda^{(0)})$. \vspace{3pt}
\State Calculate $\bP^{-1}_{\text{VV}}(r_0) = \sigma_1\iS^{\text{VtD}}_1(r_0) + \sigma_2\iS^{\text{VtD}}_2(r_0)$. \vspace{3pt}
\State Set $\tilde{r}_0 := \bP^{-1}_{\text{VV}}(r_0)$ and $q_0: = \tilde{r}_0$. \vspace{3pt}
\For{$k= 1, \cdot, N_{\text{max}}}$: \Comment{{\color{gray}\parbox[t]{.4\linewidth}{Start GMRES iterations.}}}
\State \parbox[t]{\dimexpr\linewidth-\algorithmicindent}{Generate $\lambda^{(k)}$ as the solution to the least square problem: \vspace{-0.2cm}
$$
\min\limits_{\mu \in R_k}\|\bP^{-1}_{\text{VV}}(\chi_{\iP} - \iS_{\iP}(\mu))\|_{L^{2}}, \vspace{-0.2cm}
$$
where $R_k := \lambda^{(0)} + \text{span}(q_0, q_1, \cdots, q_{k-1})$.} \vspace{1pt}
\State Set $\tilde{r}_k := \bP^{-1}_{\text{VV}}\left(\chi_{\iP}-\iS_{\iP}(\lambda^{(k)})\right)$. \vspace{3pt}
\If{$\| \tilde{r}_k \|/\| \tilde{r}_0 \| \leq \epsilon$}
\State stop the iteration, return $ \lambda = \lambda^{(k)}$.
\EndIf
\State \parbox[t]{\dimexpr\linewidth-\algorithmicindent}{Compute $ q_k :=  \bP^{-1}_{\text{VV}}(\iS_{\iP}(q_{k-1}))$ as in Steps $2$ and $4$.}
\EndFor
\alglinenoNew{alg5}
\alglinenoPush{alg5}
\end{algorithmic} 
\end{algorithm}
\begin{remark} \Rv{By definition, $ q_0=\tilde{r}_0$,  and for $k=1, \hdots, N_{\text{max}},$
$$q_k=\bP^{-1}_{\text{VV}}\iS_{\iP}(q_{k-1})=\bP^{-1}_{\text{VV}}\iS_{\iP}\left((\bP^{-1}_{\text{VV}}\iS_{\iP})^{k-1}(q_0)\right)=(\bP^{-1}_{\text{VV}}\iS_{\iP})^{k}(q_0).$$
Thus, the space $R_{k}$ in Step $8$ is the Krylov subspace corresponding to $\bP^{-1}_{\text{VV}}\iS_{\iP}$:
\begin{align*}
R_{k}=\lambda^{(0)} + \text{span}\left( \tilde{r}_0, (\bP^{-1}_{\text{VV}}\iS_{\iP})(\tilde{r}_0), \cdots, (\bP^{-1}_{\text{VV}}\iS_{\iP})^{k-1}(\tilde{r}_0)\right).  
\end{align*}}\vspace{-0.4cm}
\end{remark} 
\section{Global-in-time dual Schur (GTD-Schur) method}\label{sec4}
The dual Schur method is obtained by imposing Neumann boundary conditions for the subdomain problems,  instead of Dirichlet conditions as in the primal Schur approach.  Due to the presence of a high permeability fracture in the medium, the normal flux may not be continuous across the fracture-interface.  Thus, we introduce two variables \vspace{-0.2cm}
$$\varphi_{i}:= \bu_i\cdot\bn_{i\vert\gamma}, \; i=1, 2, \vspace{-0.2cm}$$ 
representing the normal flux from each subdomain along the fracture.
To formulate the interface problem for GTD-Schur with two unknowns $\varphi_1$ and $\varphi_2$,  we define the Neumann-to-Dirichlet operator: \vspace{-0.2cm}
\begin{align*}
\begin{array}{ccc}
\iS^{\text{NtD}}_i: L^2\left(0, T; L^2(\gamma)\right) \times L^2\left(0, T; L^2(\Omega_i)\right) \times H^{1}_{*, \gamma}(\Omega_i) & \longrightarrow &H^{1}\left(0, T; L^2(\gamma)\right), \\
(\varphi_{i}, q_{i}, p_{0, i}) & \longmapsto& (p_i)_{\vert\gamma},
\end{array} \vspace{-0.2cm}
\end{align*}
where $\left(p_i, \bu_i\right), \; i=1,2$ is the solution to the subdomain problem with Neumann conditions: \vspace{-0.2cm}
\begin{equation}
\label{subdomain_solver_DualSchur}
\left.\begin{array}{rcll}
\phi_i\partial_t{p_i}+\text{div }\bu_i&=&q_{i} &\text{ in } \Omega_i\times (0, T), \\
\bu_i&=&-\bK_i\nabla{p_i} &\text{ in } \Omega_i\times (0, T), \\
p_i&=&0 &\text{ on } \left(\partial\Omega_i \cap \partial\Omega\right) \times (0, T), \\
\bu_i \cdot \bn_i&=&\varphi_{i} &\text{ on } \gamma \times (0, T), \\
p_i(\cdot, 0)&=&p_{0, i} &\text{ in } \Omega_i. 
\end{array}\right. \vspace{-0.2cm}
\end{equation}
Next we denote by $S_{\gamma}$ the local operator on the fracture: \vspace{-0.2cm}
\begin{align*}
\begin{array}{ccc}
\iS_{\gamma}: \left (L^2\left(0, T; L^2(\gamma)\right)\right )^{2} \times L^2\left(0, T; L^2(\gamma)\right) \times H^{1}_{0}(\gamma) & \longrightarrow & H^1\left(0, T; L^2(\gamma)\right), \\
(\varphi_1, \varphi_2, q_{\gamma}, p_{0, \gamma}) & \longmapsto & p_{\gamma},
\end{array} \vspace{-0.2cm}
\end{align*}
where $\left(p_\gamma, \bu_{\gamma}\right)$ is the solution to the $(d-1)$-dimensional fracture problem: \vspace{-0.2cm}
\begin{equation}
\label{fracture_solver_DualSchur}
\begin{array}{rcll}
\phi_{\gamma}\partial_t{p_{\gamma}}+\text{div}_{\tau}\bu_{\gamma}&=& q_{\gamma} + \sum\limits^{2}_{i=1}\varphi_i & \text{ in } \gamma \times (0, T), \\
\bu_{\gamma} &=&-\bK_{\gamma}\delta\nabla_{\tau}p_{\gamma} & \text{ in } \gamma \times (0, T), \\
p_{\gamma}&=&0 &\text{ on } \partial\gamma \times (0, T), \\
p_{\gamma}(\cdot, 0)&=&p_{0, \gamma} & \text{ in } \gamma.
\end{array} \vspace{-0.2cm}
\end{equation}
The space-time interface problem is obtained by enforcing the continuity of the pressure across the fracture and is given by \vspace{-0.2cm}
\begin{equation}
\label{interface_dual_Schur}
\begin{array}{llll}
\iS_{\gamma}(\varphi_1, \varphi_2, q_{\gamma}, p_{0, \gamma}) & = & \iS^{\text{NtD}}_{1}(\varphi_{1}, q_{1}, p_{0, 1}), \text{ in } \gamma \times (0, T),  \vspace{2pt}\\
\iS_{\gamma}(\varphi_1, \varphi_2, q_{\gamma}, p_{0, \gamma}) & = & \iS^{\text{NtD}}_{2}(\varphi_{2}, q_{2}, p_{0, 2}), \text{ in } \gamma \times (0, T),
\end{array} \vspace{-0.2cm}
\end{equation}
or in compact form, \vspace{-0.2cm}
\begin{equation}
\label{interface_dual_Schur_compact_form}
\begin{array}{llll}
\iS_{\iD}\begin{pmatrix}
\varphi_1 \\
\varphi_2
\end{pmatrix} &= &\chi_{\iD}, & \; \text{in} \; \gamma \times (0, T), 
\end{array} \vspace{-0.2cm}
\end{equation}
where \vspace{-0.2cm}
\begin{equation}
\label{left_hand_DualSchur}
\iS_{\iD}\begin{pmatrix}
\varphi_1 \\
\varphi_2
\end{pmatrix} = \left(\begin{array}{c}
\iS_{\gamma}(\varphi_1, \varphi_2, 0, 0) - \iS^{\text{NtD}}_{1}(\varphi_{1}, 0, 0)  \vspace{2pt}\\
\iS_{\gamma}(\varphi_1, \varphi_2, 0, 0) -\iS^{\text{NtD}}_{2}(\varphi_{2}, 0, 0)
\end{array}\right), 
\end{equation}
and  
\begin{equation}
\label{right_hand_DualSchur}
\chi_{\iD}= \left(\begin{array}{c}
\iS^{\text{NtD}}_{1}(0, q_{1}, p_{0, 1}) -\iS_{\gamma}(0, 0, q_{\gamma}, p_{0, \gamma})  \vspace{2pt}\\
\iS^{\text{NtD}}_{2}(0, q_{2}, p_{0, 2}) -\iS_{\gamma}(0, 0, q_{\gamma}, p_{0, \gamma})
\end{array}\right). \vspace{0.1cm}
\end{equation}
The interface problem \eqref{interface_dual_Schur_compact_form} is solved iteratively, and we propose the following Dirichlet-Dirichlet preconditioner, $\bP^{-1}_{\text{DD}}$,  to enhance its convergence (cf.~Section~\ref{sec7}): \vspace{-0.2cm}
\begin{equation}
\label{interface_dual_Schur_Preconditioner}
\bP^{-1}_{\text{DD}}\left(\iS_{\iD}\begin{pmatrix}
\varphi_1 \\
\varphi_2
\end{pmatrix}\right) = {\bP^{-1}_{\text{DD}}}\left(\chi_{\iD}\right), \; \; \text{in} \; \gamma \times (0, T), \vspace{-0.2cm}
\end{equation}
where 
\begin{equation}
\label{D_D_Preconditioner}
\bP^{-1}_{\text{DD}}\begin{pmatrix}
\lambda_1 \\
\lambda_2 
\end{pmatrix} = \begin{pmatrix}
\widetilde{\iS}^{\text{DtN}}_1\left(\lambda_1\right) \vspace{2pt}\\
\widetilde{\iS}^{\text{DtN}}_2\left(\lambda_2\right)
\end{pmatrix},
\end{equation}
and $\widetilde{\iS}^{\text{DtN}}_i, \; i=1, 2$ is a Dirichlet-to-Neumann operator defined as \vspace{-0.2cm}
\begin{equation}
\label{DtN_DualSchur}
\widetilde{\iS}^{\text{DtN}}_{i}(\lambda_i) := \iS^{\text{DtN}}_i\left(\lambda_i, 0, 0\right) =  \bu_i\cdot\bn_{i\vert\gamma}. \vspace{-0.2cm}
\end{equation}
\Rv{The GTD-Schur method with the Dirichlet-Dirichlet preconditioner is outlined in Algorithm~\ref{GTD_Schur}.  The case without preconditioner follows the same steps with $\bP^{-1}_{\text{DD}}$ being replaced by the identity operator.  }
%
%
\begin{algorithm}[!ht]
\caption{GTD-Schur method with Dirichlet-Dirichlet preconditioner}
\label{GTD_Schur}
\textbf{Input:} initial guess $(\varphi^{(0)}_1, \varphi^{(0)}_2)$, stopping tolerance $ 0 <\epsilon \ll 1$, maximum number of iterations $N_{\text{max}}$. \\
\textbf{Output:} pair of space-time fracture normal fluxes $ (\varphi_1, \varphi_2).$ 
\begin{algorithmic}[1]
\State Compute 
$
\chi_{\iD}= \left(\begin{array}{c}
\iS^{\text{NtD}}_{1}(0, q_{1}, p_{0, 1}) -\iS_{\gamma}(0, 0, q_{\gamma}, p_{0, \gamma})  \vspace{2pt}\\
\iS^{\text{NtD}}_{2}(0, q_{2}, p_{0, 2}) -\iS_{\gamma}(0, 0, q_{\gamma}, p_{0, \gamma})
\end{array}\right).
$ \vspace{3pt}
\State Evaluate $\iS_{\iD}(\varphi^{(0)}_1, \varphi^{(0)}_2)=\left(\begin{array}{c}
\iS_{\gamma}(\varphi^{(0)}_1, \varphi^{(0)}_2, 0, 0) - \iS^{\text{NtD}}_{1}(\varphi^{(0)}_{1}, 0, 0)  \vspace{2pt}\\
\iS_{\gamma}(\varphi^{(0)}_1, \varphi_2^{(0)}, 0, 0) -\iS^{\text{NtD}}_{2}(\varphi^{(0)}_{2}, 0, 0)
\end{array}\right).$   \vspace{3pt}
\State Set $r_0 =(r_{0, 1}, r_{0, 2}):= \chi_{\iD} - \iS_{\iD}(\varphi^{(0)}_1, \varphi^{(0)}_2).$  \vspace{3pt}
\State Compute $\bP^{-1}_{\text{DD}}(r_0)=\begin{pmatrix}
\widetilde{\iS}^{\text{DtN}}_1\left(r_{0, 1}\right) \vspace{2pt}\\
\widetilde{\iS}^{\text{DtN}}_2\left(r_{0, 2}\right)
\end{pmatrix}.$ 
\State Set $\tilde{r}_0 = \bP^{-1}_{\text{DD}}(r_0)$ and $q_0 = \tilde{r}_0$.  \vspace{3pt}
\For{$k= 1, \cdots, N_{\text{max}}}$: \Comment{{\color{gray}\parbox[t]{.4\linewidth}{Start GMRES iterations.}}}
\State \parbox[t]{\dimexpr\linewidth-\algorithmicindent}{Generate $(\varphi^{(k)}_1, \varphi^{(k)}_2)$ as the solution to the least square problem: \vspace{-0.2cm}
$$
\min\limits_{(\psi_1, \psi_2) \in R_k}\|\bP^{-1}_{\text{DD}}(\chi_{\iD} - \iS_{\iD}(\psi_1, \psi_2))\|_{L^{2}}, \vspace{-0.2cm}
$$
where $R_k := (\varphi^{(0)}_1, \varphi^{(0)}_2) + \text{span}(q_0, q_1, \cdots, q_{k-1})$.} \vspace{1pt}
\State Set $\tilde{r}_k = \bP^{-1}_{\text{DD}}\left(\chi_{\iD}-\iS_{\iD}(\varphi^{(k)}_1, \varphi^{(k)}_2)\right)$. \vspace{3pt}
\If{$\| \tilde{r}_k \| / \| \tilde{r}_0 \| \leq \epsilon$}
\State stop the iteration, return $ (\varphi_1, \varphi_2) = (\varphi^{(k)}_1, \varphi^{(k)}_2)$.
\EndIf
\State \parbox[t]{\dimexpr\linewidth-\algorithmicindent}{Compute $q_k = \bP^{-1}_{\text{DD}}\left(\iS_{\iD}(q_{k-1})\right)$ as in Steps $2$ and $4$.} 
\EndFor
\end{algorithmic}
\end{algorithm}

\section{Global-in-time fracture-based Schur (GTF-Schur) method} \label{sec5}
The primal and dual Schur methods generally require suitable preconditioners to achieve satisfactory convergence speed.  Though the number of iterations is reduced with preconditioning,  additional subdomain problems need to be solved.  It would be desirable to develop a DD method that converges fast without any preconditioners.  By combining the ideas of GTP-Schur and GTD-Schur, we derive the GTF-Schur method whose space-time interface operator is closed to identity operator, thus,  making the new interface problem better-conditioned.
Instead of having two interface unknowns as in the GTD- Schur method, only one term $\varphi := \sum\limits^{2}_{i=1}\bu_i\cdot\bn_{i\vert\gamma}$ representing the jump of the normal flux across the fracture will be introduced. The fracture pressure $p_{\gamma}$ is then recovered by solving the fracture problem \eqref{interface_problem_Schur_b} provided the new unknown. Toward this end, we define the solution operator \vspace{-0.2cm}
\begin{align*}
\begin{array}{ccc}
\widehat{\iS}_{\gamma}:  L^2\left(0, T; L^2(\gamma)\right) \times L^2\left(0, T; L^2(\gamma)\right) \times H^{1}_{0}(\gamma) & \longrightarrow &H^1\left(0, T; L^2(\gamma)\right), \\
(\varphi, q_{\gamma}, p_{0, \gamma}) & \longmapsto & p_{\gamma},
\end{array} \vspace{-0.2cm}
\end{align*}
\noindent
where $\left(p_{\gamma}, \bu_{\gamma}\right)$ is the solution to the flow problem on the fracture: \vspace{-0.2cm}
\begin{equation}
\label{fracture_solver_ModifiedSchur}
\begin{array}{rcll}
\phi_{\gamma}\partial_t{p_{\gamma}}+\text{div}_{\tau }\bu_{\gamma}&=& q_{\gamma} + \varphi & \text{ in } \gamma \times (0, T), \\
\bu_{\gamma} &=&-\bK_{\gamma}\delta\nabla_{\tau}p_{\gamma} & \text{ in } \gamma \times (0, T), \\
p_{\gamma}&=&0 &\text{ on } \partial\gamma \times (0, T), \\
p_{\gamma}(\cdot, 0)&=&p_{0, \gamma} & \text{ in } \gamma. 
\end{array} \vspace{-0.2cm}
\end{equation}
Using $p_{\gamma}=\widehat{\iS}_{\gamma}(\varphi, q_{\gamma}, p_{0, \gamma})$ as Dirichlet boundary data on the fracture-interface,  \vspace{-0.2cm}
$$ p_i=p_{\gamma}, \quad \text{on} \; \gamma \times (0, T),  \vspace{0.1cm}
$$
we solve the subdomain problem~ \eqref{subdom_solver_Schur} to obtain $(p_{i}, \bu_{i})$,  from which the normal flux is computed: \vspace{-0.2cm}
\begin{align*}
\iS^{\text{DtN}}_i({\widehat{\iS}_{\gamma}}(\varphi, q_{\gamma}, p_{0, \gamma}), q_{i}, p_{0, i}) = \bu_i\cdot\bn_{i\vert\gamma}, \; i =1, 2, \vspace{-0.2cm}
\end{align*}
where $\iS^{\text{DtN}}$ is the same Dirichlet-to-Neumann operator as in GTP-Schur.
Finally, the interface problem for GTF-Schur is obtained by matching $\varphi$ with the total normal fluxes: \vspace{-0.2cm}
\begin{equation} \label{interface_problem_ModifiedSchur}
\begin{array}{lll}
\varphi & = & \sum\limits^{2}_{i =1} \iS^{\text{DtN}}_i({\widehat{\iS}_{\gamma}}(\varphi, q_{\gamma}, p_{0, \gamma}), q_{i}, p_{0, i}), \; \text{in} \; \gamma \times (0, T),
\end{array} \vspace{-0.2cm}
\end{equation}
or in compact form, \vspace{-0.2cm}
\begin{equation}
\label{IP_FS_compact_form}
\begin{array}{llll}
\iS_{\iF}\left (\varphi \right ) &= &\chi_{\iF}, & \; \text{in} \; \gamma \times (0, T), 
\end{array} \vspace{-0.2cm}
\end{equation}
where $$\iS_{\iF}\left (\varphi \right ) =\sum\limits^{2}_{i =1}  \iS^{\text{DtN}}_i({\widehat{\iS}_{\gamma}}(\varphi, 0, 0), 0, 0), \;\chi_{\iF}= \sum\limits^{2}_{i =1}  \iS^{\text{DtN}}_i({\widehat{\iS}_{\gamma}}(0, q_{\gamma}, p_{0, \gamma}), q_{i}, p_{0, i}).$$ 
Again, we solve the interface problem~\eqref{IP_FS_compact_form} iteratively using GMRES \Rv{(without any preconditioner) as summarized in Algorithm~\ref{GTF_Schur}}.  Numerical performance of GTF-Schur will be discussed and compared with GTP-Schur and GTD-Schur in Section~\ref{sec7}.
%
%
\begin{algorithm}[!ht]
\caption{GTF-Schur method}
\label{GTF_Schur}
\textbf{Input:} initial guess $\varphi^{(0)}$, stopping tolerance $ 0 <\epsilon \ll 1$, maximum number of iterations $N_{\text{max}}$. \\
\textbf{Output:} space-time total normal flux $\varphi$.  
\begin{algorithmic}[1]
\State Compute 
$\chi_{\iF}= \sum\limits^{2}_{i =1}  \iS^{\text{DtN}}_i({\widehat{\iS}_{\gamma}}(0, q_{\gamma}, p_{0, \gamma}), q_{i}, p_{0, i}).$
\State Evaluate  $\iS_{\iF}(\varphi^{(0)}) =\sum\limits^{2}_{i =1}  \iS^{\text{DtN}}_i({\widehat{\iS}_{\gamma}}(\varphi^{(0)}, 0, 0), 0, 0).$  \vspace{3pt}
\State Set $r_0 = \chi_{\iF} - \iS_{\iF}(\varphi^{(0)})$. \vspace{3pt}
\For{$k= 1, \cdots, N_{\text{max}}}$: \Comment{{\color{gray}\parbox[t]{.4\linewidth}{Start GMRES iterations.}}}
\State \parbox[t]{\dimexpr\linewidth-\algorithmicindent}{Generate $\varphi^{(k)}$ as a solution to the least square problem: \vspace{-0.2cm}
$$
\min\limits_{\psi \in R_k}\|\chi_{\iF} - \iS_{\iF}(\psi)\|_{L^{2}}, \vspace{-0.2cm}
$$
where $R_k := \varphi^{(0)} + \text{span}(r_0, \iS_{\iF}(r_0), \cdots, \iS^{k-1}_{\iF}(r_0))$.} \vspace{1pt}
\State Set $r_k = \chi_{\iF}-\iS_{\iF}(\varphi^{(k)})$. \vspace{3pt}
\If{$\| r_k \| / \| r_0 \| \leq \epsilon$}
\State stop the iteration, return $ \varphi = \varphi^{(k)}$.
\EndIf
\State Compute $\iS^{k}_{\iF}(r_0) = \iS_{\iF}\left(\iS^{k-1}_{\iF}(r_0)\right)$ as in Step $2$.
\EndFor
\end{algorithmic}
\end{algorithm} \vspace{-0.4cm}
\section{Nonconforming discretization in time}\label{sec6}
All three DD methods presented in previous sections are globally in time, i.e., the subdomain problems are solved over the whole time interval at each iteration and space-time information is exchanged on the fracture-interface.  Thus it is possible to use different time steps in the fracture and in the rock matrix.  In this section, we derive the semidiscrete interface problem for the proposed DD methods with nonconforming time grids.

Let $\mathcal{T}_1, \mathcal{T}_2,$ and $\mathcal{T}_{\gamma}$ be three different partitions of the time interval $(0, T]$ into subintervals $J^{i}_{m} = \left(t^{i}_{m-1}, t^{i}_m\right]$ for $m=1, \cdots, M_i,$ and $i=1, 2, \gamma$ (see Figure~\ref{nonconforming_time}). For simplicity,  we consider uniform partitions and denote by $\Delta{t}_i, \; i=1, 2, \gamma$,  the corresponding time steps such that $\Delta{t}_{\gamma} \ll \Delta{t}_i, \; i=1,2$ (note that the fracture is assumed to have much larger permeability than the surround domain).  We use the backward Euler method to discretize the problem in time.  The same idea can be generalized to higher order methods \cite{23}.

We denote by $P_0\left(\iT_i, \; L^2(\gamma)\right)$ the space of functions which are piecewise constant in time on grid $\iT_i$ with values in $L^2(\gamma)$: \vspace{-0.2cm}
\begin{align*}
P_0\left(\iT_i, \; L^2(\gamma)\right) = \left\{ \psi: (0, T) \rightarrow L^2(\gamma), \psi \; \text{is} \; \text{constant on} \; J,  \; \forall J \in \iT_i\right\}. \vspace{-0.2cm}
\end{align*} 
In order to exchange data on the space-time interface between different time grids $\iT_{i}$ and $\iT_{j}$ (for $i, j$ in $\{1, 2, \gamma\}$),  we use the $L^2$ projection $\Pi_{ji}$ from $P_0\left(\iT_i, \; L^2(\gamma)\right)$ to $P_0\left(\iT_j, \; L^2(\gamma)\right)$: for $\psi \in P_0\left(\iT_i, \; L^2(\gamma)\right),$ $\Pi_{ji}\psi_{\vert{J^{j}_m}}$ is the average value of $\psi$ on $J^{j}_{m},$ for $m= 1, \cdots, M_j$. 
\begin{figure}[t!]
\centering
\includegraphics[scale=0.65]{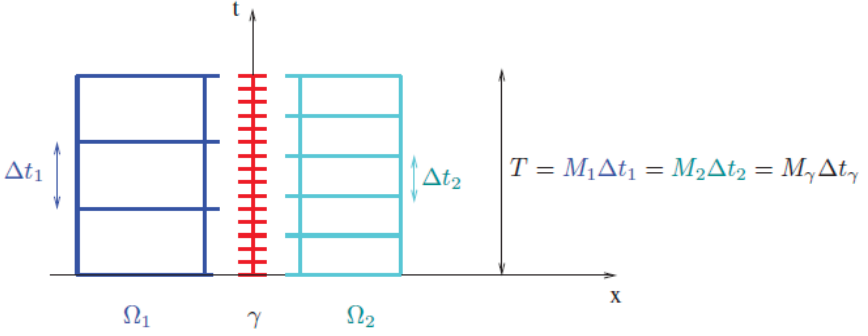}
\caption{Nonconforming time grids in the rock matrix and in the fracture.}
\label{nonconforming_time} \vspace{-0.3cm}
\end{figure}
\subsection{GTP-Schur method} \label{subsec6a}
The unknown $\lambda$ in \eqref{interface_problem_Schur_a} is piecewise constant in time on grid $\mathcal{T}_{\gamma}$ as it represents the pressure on the fracture.  In order to obtain Dirichlet boundary data for the subdomain problem \eqref{subdom_solver_Schur}, we project $\lambda$ into $P_0\left(\mathcal{T}_i, \; L^2(\gamma)\right)$: 
$ p_i = \Pi_{i\gamma}\left(\lambda\right) \; \; \text{on}, \; i= 1, 2.  $
The semidiscrete counterpart of the interface problem \eqref{interface_problem_Schur_a} is obtained by weakly enforcing the fracture problem over each time subinterval of $\mathcal{T}_{\gamma}$ as follows: \vspace{-0.2cm}
\begin{equation}
\begin{array}{rl}
\phi_{\gamma}\left(\lambda^{m+1}-\lambda^m\right) + {\int}^{t^{m+1}_{\gamma}}_{t^{m}_{\gamma}}\text{div}_{\tau}\; \bu^{m+1}_{\gamma} &= {\int}^{t^{m+1}_{\gamma}}_{t^{m}_{\gamma}}\left(\sum\limits^2_{i=1}\Pi_{\gamma{i}}\left(\iS^{\text{DtN}}_i\left(\Pi_{i\gamma}(\lambda), q_{i}, p_{0, i}\right)\right)\right),\vspace{2pt}\\
\bu^{m+1}_{\gamma} &= -\bK_{\gamma}\delta{\nabla}_{\tau}\lambda^{m+1},
\end{array} \label{eq:semiGTPS}
\end{equation}
in $\gamma$,  for $m=0, \cdots, M_{\gamma}-1$.  Problem~\eqref{eq:semiGTPS} is completed with the initial and boundary conditions: $ \lambda^0 = p_{0, \gamma},  \; \text{in} \; \gamma$ and $\lambda^{m+1} = 0, \; \text{on} \; \partial\gamma, $ for $m=0, \cdots, M_{\gamma}-1$.

To compute the semidiscrete Ventcel-Ventcel preconditioner, which is still denoted by $P^{-1}_{\text{VV}}$, we first project the data $\theta \in P_0\left(\mathcal{T}_{\gamma}, L^2\left(\gamma\right)\right)$ onto the subdomain grid $\mathcal{T}_{i}, \; i=1, 2$ to solve the subdomain problem with Ventcel conditions~\eqref{subdom_solver_Schur_preconditioner}. 
Then $P^{-1}_{\text{VV}}$ is obtained by projecting the trace of the subdomain pressure on the fracture-interface from $\mathcal{T}_i$ onto $\mathcal{T}_{\gamma}$: \vspace{-0.2cm}
\begin{equation}
P^{-1}_{\text{VV}}\left(\theta\right) := \sum\limits^{2}_{i=1}\sigma_i\Pi_{\gamma{i}}\left(\iS^{\text{VtD}}_i\left(\Pi_{i\gamma}(\theta)\right)\right). \vspace{-0.2cm}
\end{equation}
%
\subsection{GTD-Schur method} \label{subsec6b}
The two interface unknowns $\varphi_1$ and $\varphi_2$ are piecewise constant in time on the fine grid $\mathcal{T}_{\gamma}$: $\varphi_i\in P_0\left(\mathcal{T}_{\gamma}, L^2\left(\gamma\right)\right)$ for $i=1,2$.  In order to obtain Neumann boundary data for the subdomain problem \eqref{subdomain_solver_DualSchur}, we project $\varphi_i$ into $P_0\left(\mathcal{T}_i, \; L^2(\gamma)\right)$:
$ \bu_i \cdot \bn_{i}= \Pi_{i\gamma}\left(\varphi_{i}\right) \; \; \text{on}, \; i= 1, 2.  $
The semidiscrete counterpart of the interface problem \eqref{interface_dual_Schur} is defined on $\mathcal{T}_{\gamma}$ as follows:
\vspace{-0.2cm}
\begin{equation}
\begin{array}{lr}
{\int}^{t^{m+1}_{\gamma}}_{t^{m}_{\gamma}} \iS_{\gamma}\left(\varphi_1, \varphi_2, q_{\gamma}, p_{0, \gamma}\right)&= {\int}^{t^{m+1}_{\gamma}}_{t^{m}_{\gamma}}\Pi_{\gamma{1}}\left(\iS^{\text{NtD}}_1\left(\Pi_{1\gamma}(\varphi_1\right), q_{1}, p_{0, 1}\right),\vspace{0.15cm} \\ 
{\int}^{t^{m+1}_{\gamma}}_{t^{m}_{\gamma}} \iS_{\gamma}\left(\varphi_1, \varphi_2, q_{\gamma}, p_{0, \gamma}\right) &= {\int}^{t^{m+1}_{\gamma}}_{t^{m}_{\gamma}} \Pi_{\gamma{2}}\left(\iS^{\text{NtD}}_2\left(\Pi_{2\gamma}(\varphi_2\right), q_{2}, p_{0, 2}\right), 
\end{array}    \vspace{-0.2cm}
\end{equation}
in $\gamma$, for $m=0, \cdots, M_{\gamma}-1$.\vspace{0.1cm}

The semidiscrete Dirichlet-Dirichlet preconditioner $\bP_{DD}^{-1}$ is computed by \vspace{-0.2cm}
\begin{equation}
\bP_{DD}^{-1}\begin{pmatrix}
\lambda_1 \\
\lambda_2
\end{pmatrix} = \left(\begin{array}{c}
\Pi_{\gamma{1}}\left(\widetilde{\iS}^{\text{DtN}}_1\left(\Pi_{1\gamma}(\lambda_1)\right)\right) \vspace{2pt}\\
\Pi_{\gamma{2}}\left(\widetilde{\iS}^{\text{DtN}}_2\left(\Pi_{2\gamma}(\lambda_2)\right)\right)
\end{array}\right),  \vspace{-0.2cm}
\end{equation}
\noindent
in which we first solve the subdomain problems with Dirichlet data projected from the $\mathcal{T}_{\gamma}$ onto $\mathcal{T}_{i}, \; i=1,2$, then extract the normal flux along the fracture and project backward from $\mathcal{T}_{i}$ onto $\mathcal{T}_{\gamma}$.

\subsection{GTF-Schur method} \label{subsec6c}
The interface unknown $\varphi$ in this case represents the total normal flux, and again,  it is piecewise constant in time on $\mathcal{T}_{\gamma}$: $\varphi \in P_0\left(\mathcal{T}_{\gamma}, L^2\left(\gamma\right)\right)$.  Solving the fracture problem \eqref{fracture_solver_ModifiedSchur} with $\varphi$, we obtain $p_{\gamma} = \widehat{\iS}_{\gamma}\left(\varphi, q_{\gamma}, p_{0, \gamma}\right) \in P_0\left(\mathcal{T}_{\gamma}, L^2\left(\gamma\right)\right)$.  As for GTP-Schur,  the fracture pressure $p_{\gamma}$ is projected to $\iT_{i}$, for $i=1,2,$ to give Dirichlet data for solving the subdomain problems.  The semidiscrete counterpart of \eqref{interface_problem_ModifiedSchur} is then defined on $\mathcal{T}_{\gamma}$ as follows:
%
\vspace{-0.2cm}
\begin{equation}
\begin{array}{ll}
\varphi^{m} - &\sum\limits^{2}_{i=1}\Pi_{\gamma{i}}\left(\iS^{\text{DtN}}_i\left(\Pi_{i\gamma}\left({\widehat{\iS}_{\gamma}}(\varphi, 0, 0)\right), 0, 0\right)\right)\vert_{J^{\gamma}_{m}} \\
& \hspace{2cm}= \sum\limits^{2}_{i=1}\Pi_{\gamma{i}}\left(\iS^{\text{DtN}}_i\left(0, q_{\gamma}, p_{0, \gamma}, q_{i}, p_{0, i}\right)\right)\vert_{J^{\gamma}_m},
\end{array}  \vspace{-0.2cm}
\end{equation}
on $\gamma$,  for $m=0, \cdots, M_{\gamma}-1$.
%
\section{Numerical results}\label{sec7} 
We study and compare the convergence and accuracy in time of four global-in-time DD methods: GTP-Schur with Neumann-Neumann (N-N) or Ventcel-Ventcel (V-V) preconditioners,  GTD-Schur with Dirichlet-Dirichlet (D-D) preconditioner,  GTF-Schur,  and GTO-Schwarz.  We refer to~\cite{26} for the detailed derivation and formulation of the GTO-Schwarz method and optimized parameters.  

Two test cases are considered: Test case 1 with a non-immersed fracture (i.e., the fracture cuts through the rock matrix) and Test case 2 with a partially immersed fracture.  For both cases,  we assume that $\bK_i = \mathsf{k}_iI,$ for $i=1, 2, f,$ where $\mathsf{k}_1 = \mathsf{k}_2 =1$ and $\mathsf{k}_{\gamma} = 10^3$.  
For spatial discretization, we consider mixed finite elements with the lowest order Raviart–Thomas space on a uniform, conforming triangular mesh of size $h$. 
We remark that the focus of this work is local time stepping; nonconforming spatial meshes will be the topic of our future work.  The interface problem for each method is solved iteratively using GMRES with a random initial guess; the iteration is stopped when the residual error is less than $10^{-6}$ (Test case~1) or $10^{-8}$ (Test case 2).  All computed errors are relative space-time errors in the space $L^2(0, T; L^2(\mathcal{O}))$-norm, where $\mathcal{O}$ is either $\Omega_{1}$,  $\Omega_{2}$,  or $\gamma$.  To compare the convergence of the iterative algorithms (with or without preconditioners),  we count the number of subdomain solves instead of the number of iterations. Note that one iteration of GTP-Schur or GTD-Schur with a preconditioner costs twice as much as one iteration of
the respective method with no preconditioner (in terms of number of subdomain solves).
\subsection{Test case 1: with a non-immersed fracture}\label{subsec7}
\begin{figure}[h]
\centering
\vspace{-0.3cm}
\begin{minipage}{.5\textwidth}
  \centering
  \includegraphics[width=0.95\linewidth]{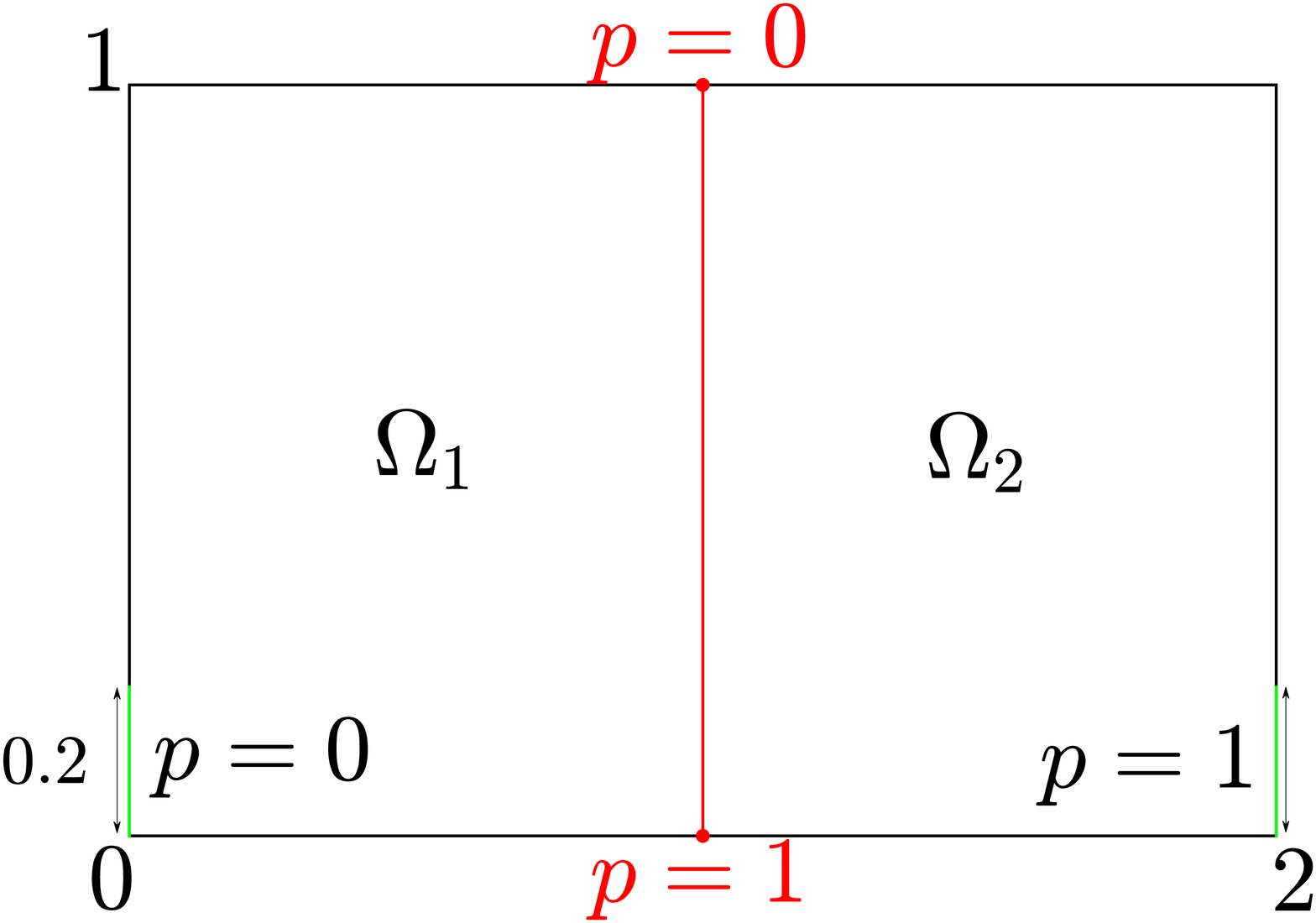}
  \label{fig:test1}
\end{minipage}%
\begin{minipage}{.5\textwidth}
  \centering
  \includegraphics[width=0.95\linewidth]{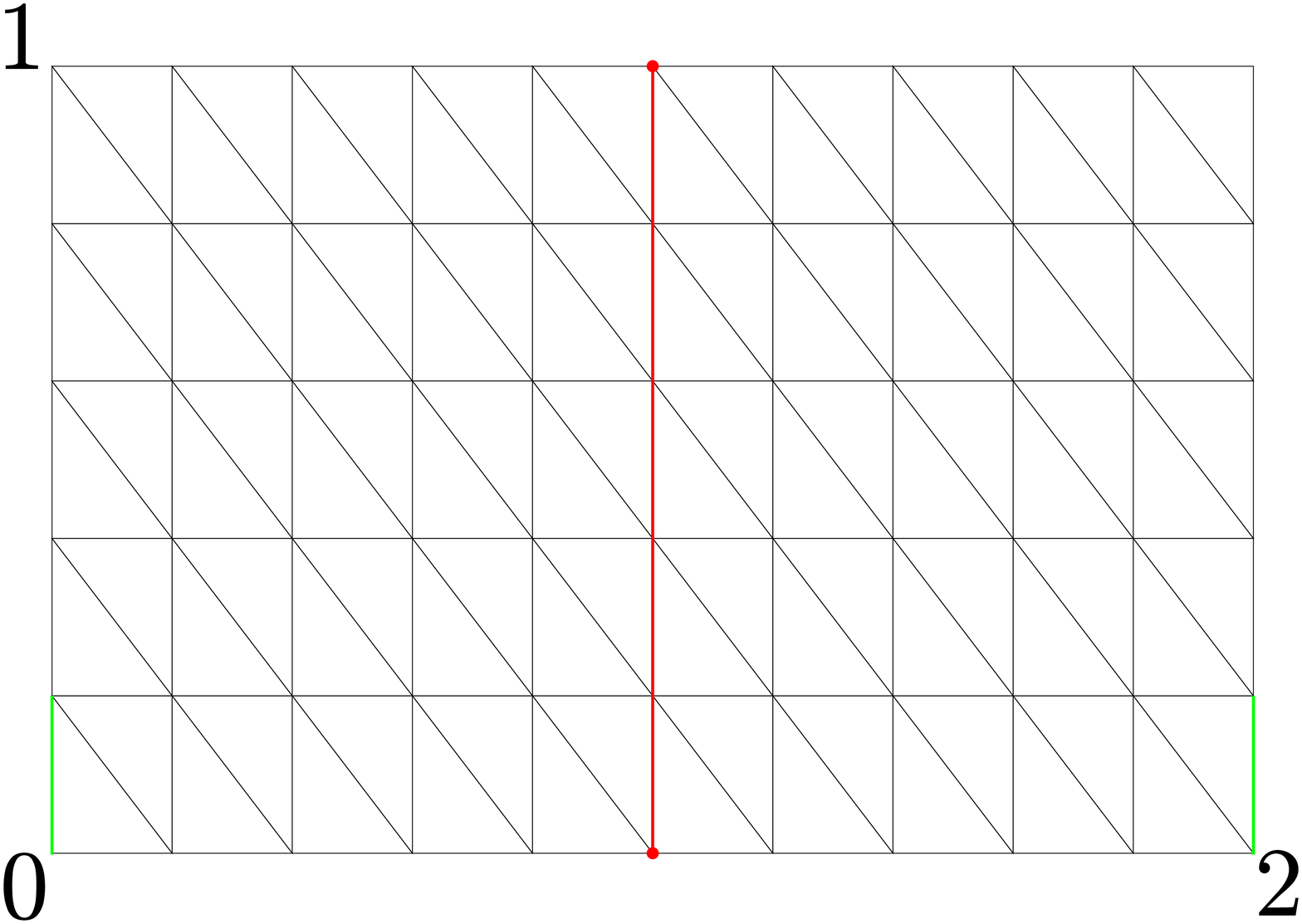}
  \label{fig:test2}
\end{minipage}
\caption{[Test case 1] (Left) Geometry and boundary conditions of the test case.  (Right) Example of an uniform triangular mesh for spatial discretization.}
\label{Non_Immersed} \vspace{-0.4cm}
\end{figure}

\indent
The domain of calculation $\Omega=(0,2) \times (0,1)$ is divided into two equally sized subdomains by a fracture of width $\delta = 0.001$ parallel to the $y$-axis (see Figure~\ref{Non_Immersed}).  For the boundary conditions, we impose $p =1$ at the bottom and $p=0$ at the top of the fracture.  On the external boundaries of the subdomains, a no flow boundary condition is imposed except on the lower fifth (length 0.2) of both lateral sides where a Dirichlet condition is imposed: $p = 1$ on the right and $p = 0$ on the left. 

We first consider the {\em conforming} time step case to verify the errors and compare the convergence of the four global-in-time DD algorithms.  We fix the final time $T=0.5$,  the spatial mesh $h=1/50$, and vary the time step sizes $\Delta t_{i}=\Delta t$ for $i=1,2,\gamma$.  The errors are computed using a reference solution on a fine time step $\Delta{t}_{\text{ref}} = T/2000$.  Table~\ref{Err_Concentration_Velocity_Non_Immersed} shows the $L^2$ errors for the pressure and velocity computed once GMRES converges.  Note that all methods produce nearly the same approximate solutions since the same time step is imposed in the fracture and in the subdomains.  From this table, first order convergence in time is observed for both pressure and velocity.
%
\begin{table}[h]
\centering
\begin{minipage}{\textwidth}
\begin{tabular*}{\textwidth}{@{\extracolsep{\fill}}lcccccc@{\extracolsep{\fill}}}
\toprule%
& \multicolumn{3}{@{}c@{}}{Errors for pressure} & \multicolumn{3}{@{}c@{}}{Errors for velocity} \\ \cmidrule{2-4}\cmidrule{5-7}%
$\Delta{t}$ & $\Omega_1$ & $\Omega_2$ & $\gamma$ & $\Omega_1$ & $\Omega_2$ & $\gamma$ \\
\midrule
T/4  & 6.76e-02 & 6.82e-02 & 3.29e-02 &4.96e-02 & 9.24e-02 & 5.47e-02\\
&&&&& \\
T/8 & 3.55e-02 & 3.57e-02   & 1.59e-02  & 2.56e-02  & 4.87e-02 & 2.64e-02 \\
&$[0.92]$ & $[0.93]$ & $[1.05]$ & $[0.95]$ & $[0.92]$ & $[1.05]$ \\
&&&&& \\
T/16 & 1.81e-02 & 1.81e-02  & 7.73e-03  & 1.30e-02  & 2.49e-02  & 1.28e-02  \\
& $[0.97]$ & $[0.98]$ & $[1.04]$ & $[0.97]$ & $[0.96]$ & $[1.04]$ \\
&&&&& \\
T/32 & 9.06e-03  & 9.07e-03  & 3.76e-03  & 6.52e-03  & 1.24e-02  & 6.24e-03  \\
& $[0.99]$ & $[0.99]$ & $[1.03]$ & $[0.99]$ & $[1.00]$ & $[1.03]$ \\
\botrule
\end{tabular*}
\caption{[Test case 1] Relative $L^2$ errors of the pressure and velocity with {\em conforming} time steps.  The corresponding convergence rates are shown in square brackets.}\label{Err_Concentration_Velocity_Non_Immersed}
\end{minipage} \vspace{-0.4cm}
\end{table}
In Table~\ref{GMRES_iterations_Non_Immersed_conforming}, we report the number of subdomain solves needed to obtain such errors.  In particular, we stop GMRES when the relative residual is smaller than $10^{-6}$.  
For GTP-Schur, we see that without preconditioner,  the convergence is extremely slow and deteriorates as the time step decreases.  With V-V preconditioner,  the number of iterations is significantly reduced and independent of the time step size.  For GTD-Schur,  even without a preconditioner, the performance is much better than that of GTP-Schur, and applying D-D preconditioner results in a comparable result as GTP-Schur with V-V preconditioner.  Importantly, GTF-Schur works remarkably well with no preconditioner needed, and in terms of computational cost (or subdomain solves), it is the only Schur type method that can compete with GTO-Schwarz. 
\begin{table}[h] 
\centering
\setlength{\extrarowheight}{4pt}
	\begin{tabular}{llcccc } \hline 
		& $\Delta t$ & T/4  & T/8  & T/16  & T/32 \\ \hline
		\multirow{3}{*}{GTP-Schur} & with no precond.  & $191$ & $282$ & $331$ & $407$ \\ 
		& with N-N precond. & $78$ & $92$  & $102$  & $140$  \\  
		& with V-V precond. & $10$&$12$& $12$& $12$  \\ \hline
		\multirow{2}{*}{GTD-Schur} & with no precond.  & $33$  & $34$  & $33$ & $33$ \\ 
		& with D-D precond.  &$16$&$16$&$16$&$16$ \\   \hline
		\multicolumn{2}{l}{GTF-Schur} &$8$&$8$&$8$& $8$  \\ \hline
		\multicolumn{2}{l}{GTO-Schwarz} &$6$&$6$&$6$&$6$  \\ \hline
	\end{tabular}  \vspace{0.2cm}
	\caption{[Test case 1] Numbers of subdomain solves when {\em conforming} time steps are used; the tolerance for GMRES is set to be $10^{-6}$. } 
\label{GMRES_iterations_Non_Immersed_conforming} \vspace{-0.4cm}
\end{table}

Next we investigate the case with {\em nonconforming} time grids.  We only consider GTP-Schur with V-V preconditioner, GTD-Schur with D-D preconditioner, GTF-Schur and GTO-Schwarz since they give fastest convergence.  The diffusion coefficients in the subdomains are the same and smaller than that in the fracture,  thus we impose the same large time step in the subdomains and a smaller one in the fracture: $\Delta{t}_1 =  \Delta{t}_2 = 4\Delta{t}_{\gamma}$. 
We show the relative errors of the pressure and velocity in Table~\ref{Err_Pressure_Non_Immersed_Nonconforming} and~\ref{Err_Velocity_Non_Immersed_Nonconforming}, respectively. We see that these methods still preserve the first order of convergence in time when we have nonconforming discretization in time. However,
due to the nonconforming time projections, the errors are different between the following two groups:
\begin{itemize}
\item Group 1: GTP-Schur with V-V preconditioner, and GTO-Schwarz, 
\item Group 2: GTD-Schur with D-D preconditioner,  and GTF-Schur.  
\end{itemize}
\begin{table}[!http]
\centering
\begin{minipage}{\textwidth}
\begin{tabular*}{\textwidth}{@{\extracolsep{\fill}}cccccccc@{\extracolsep{\fill}}}
\toprule%
& & \multicolumn{3}{@{}c@{}}{{GTP-Schur with V-V precond.}} & \multicolumn{3}{@{}c@{}}{GTD-Schur with D-D precond.}  \\ 
& & \multicolumn{3}{@{}c@{}}{{GTO-Schwarz}} & \multicolumn{3}{@{}c@{}}{GTF-Schur}  \\ \cmidrule{3-5} \cmidrule{6-8}%
$\Delta{t}_{i}$ & $\Delta{t}_{\gamma}$ & $\Omega_1$ & $\Omega_2$ & $\gamma$ & $\Omega_1$ & $\Omega_2$ & $\gamma$\\
\midrule
T/4  & T/16 & 6.76e-02 & 6.82e-02 & 3.29e-02 & 6.34e-02 & 6.62e-02 & 1.29e-02 \\
&&& \\
T/8 & T/32  & 3.55e-02 & 3.57e-02& 1.59e-02 & 3.27e-02 & 3.43e-02 & 6.25e-03\\
& &$[0.92]$ & $[0.93]$ & $[1.05]$ & $[0.95]$ & $[0.95]$ & $[1.04]$  \\
&&& \\
T/16 & T/64 &  1.81e-02 & 1.81e-02&  7.73e-03 & 1.65e-02 & 1.73e-02 & 3.01e-03\\
& & $[0.97]$ & $[0.98]$ & $[1.04]$ & $[0.98]$ & $[0.99]$ & $[1.05]$   \\
&&& \\
T/32 & T/128 &  9.06e-03 & 9.07e-03&  3.76e-03 & 8.22e-03 & 8.64e-03 & 1.42e-03\\
& & $[0.99]$ & $[0.99]$ & $[1.03]$  & $[1.00]$ & $[1.00]$ & $[1.08]$ \\
\botrule
\end{tabular*}
\caption{ [Test case 1] Relative $L^2$ errors of the {\em pressure} with nonconforming time grids. The corresponding convergence rates are shown in square brackets.}\label{Err_Pressure_Non_Immersed_Nonconforming}
\end{minipage} \vspace{-0.4cm}
\end{table}

\begin{table}[!http]
\centering
\begin{minipage}{\textwidth}
\begin{tabular*}{\textwidth}{@{\extracolsep{\fill}}cccccccc@{\extracolsep{\fill}}}
\toprule%
& & \multicolumn{3}{@{}c@{}}{{GT-Schur with V-V precond.}} & \multicolumn{3}{@{}c@{}}{GTD-Schur with D-D precond.}  \\ 
& & \multicolumn{3}{@{}c@{}}{{GTO-Schwarz}} & \multicolumn{3}{@{}c@{}}{GTF-Schur}  \\ \cmidrule{3-5} \cmidrule{6-8}%
$\Delta{t}_{i}$ & $\Delta{t}_{\gamma}$ & $\Omega_1$ & $\Omega_2$ & $\gamma$ & $\Omega_1$ & $\Omega_2$ & $\gamma$\\
\midrule
T/4  & T/16 & 4.96e-02 & 9.24e-02 & 5.47e-02 & 4.73e-02 & 9.38e-02 & 2.21e-02 \\
&&& \\
T/8 & T/32 & 2.56e-02 & 4.87e-02 & 2.64e-02 & 2.41e-02 & 4.87e-02 & 1.06e-02\\
&& $[0.95]$ & $[0.92]$ & $[1.05]$ & $[0.97]$ & $[0.95]$ & $[1.06]$ \\
&&& \\
T/16 & T/64 & 1.30e-02 & 2.49e-02&  1.28e-02 & 1.21e-02 & 2.47e-02 & 5.09e-03\\
& & $[0.97]$ & $[0.96]$ & $[1.04]$ & $[0.99]$ & $[0.98]$ & $[1.05]$\\
&&& \\
T/32 & T/128 &  6.52e-03 & 1.24e-02&  6.24e-03 & 6.05e-03 & 1.23e-02 & 2.41e-03\\
&& $[0.99]$ & $[1.00]$ & $[1.03]$ & $[1.00]$ & $[1.00]$ & $[1.08]$ \\
\botrule
\end{tabular*}
\caption{[Test case 1] Relative $L^2$ errors of the {\em velocity} with nonconforming time grids. The corresponding convergence rates are shown in square brackets.}\label{Err_Velocity_Non_Immersed_Nonconforming}
\end{minipage} \vspace{-0.4cm}
\end{table}


It can be observed by comparing with Table~\ref{Err_Concentration_Velocity_Non_Immersed} that the errors in the fracture for both pressure and velocity obtained from Group 1 follow the coarse time grid in the subdomains.  This behavior was observed numerically in~\cite{26} for the GTO-Schwarz method.  It is due to the fact that for GTO-Schwarz and GTP-Schur with V-V preconditioner, the fracture problem is treated as the Ventcel boundary condition for the subdomain problems.  Consequently,  the approximate fracture pressure follows the coarse time grid in the subdomains. However, for the methods in Group 2,  it can be seen that the errors in the fracture are smaller and are closer to that of the fine time grid. This is because we separate the fracture problem and the subdomain problems,  and the fracture problem is actually solved on the fine time grid. 

We now analyze the convergence of the four algorithms. Table~\ref{GMRES_iterations_Non_Immersed_nonconforming} shows the number of subdomain solves for each method to reach the relative residual smaller than $10^{-6}$. We can see that the obtained numbers are almost the same as those in Table~\ref{GMRES_iterations_Non_Immersed_conforming} and are not affected by the small time steps in the fracture. Hence, these methods are suitable for using nonconforming discretization in time. From the accuracy and convergence of the four methods in this test case, it appears that GTF-Schur is the most effective method which converges fast and preserves the accuracy in time in the fracture with smaller time steps. 

\begin{table}[!http]
\centering
\setlength{\extrarowheight}{4pt}
\begin{tabular*}{\textwidth}{@{\extracolsep{\fill}}lccccc@{\extracolsep{\fill}}}
\toprule%
& $\Delta{t}_1 = \Delta{t}_2$ & T/4 & T/8 & T/16 & T/32 \\
& $\Delta{t}_{\gamma}$ & T/16 & T/32 & T/64 & T/128 \\
\midrule
GTP-Schur with V-V precond. && $12$&$12$& $12$& $14$ \\
GTD-Schur with D-D precond. &&$16$&$16$&$16$&$16$ \\
GTF-Schur &&$8$&$8$&$8$& $8$ \\
GTO-Schwarz &&$6$&$6$&$6$&$6$ \\
\botrule
\end{tabular*}
\caption{[Test case 1] Numbers of subdomain solves when {\em nonconforming} time steps are used; the tolerance for GMRES is set to be $10^{-6}$. } 
\label{GMRES_iterations_Non_Immersed_nonconforming} \vspace{-0.4cm}
\end{table}
\subsection{Test case 2: with a partially immersed fracture} \label{subsec8}
\begin{figure}[h]
\centering
\vspace{-0.3cm}
\begin{minipage}{.5\textwidth}
  \centering
  \includegraphics[width=0.97\linewidth]{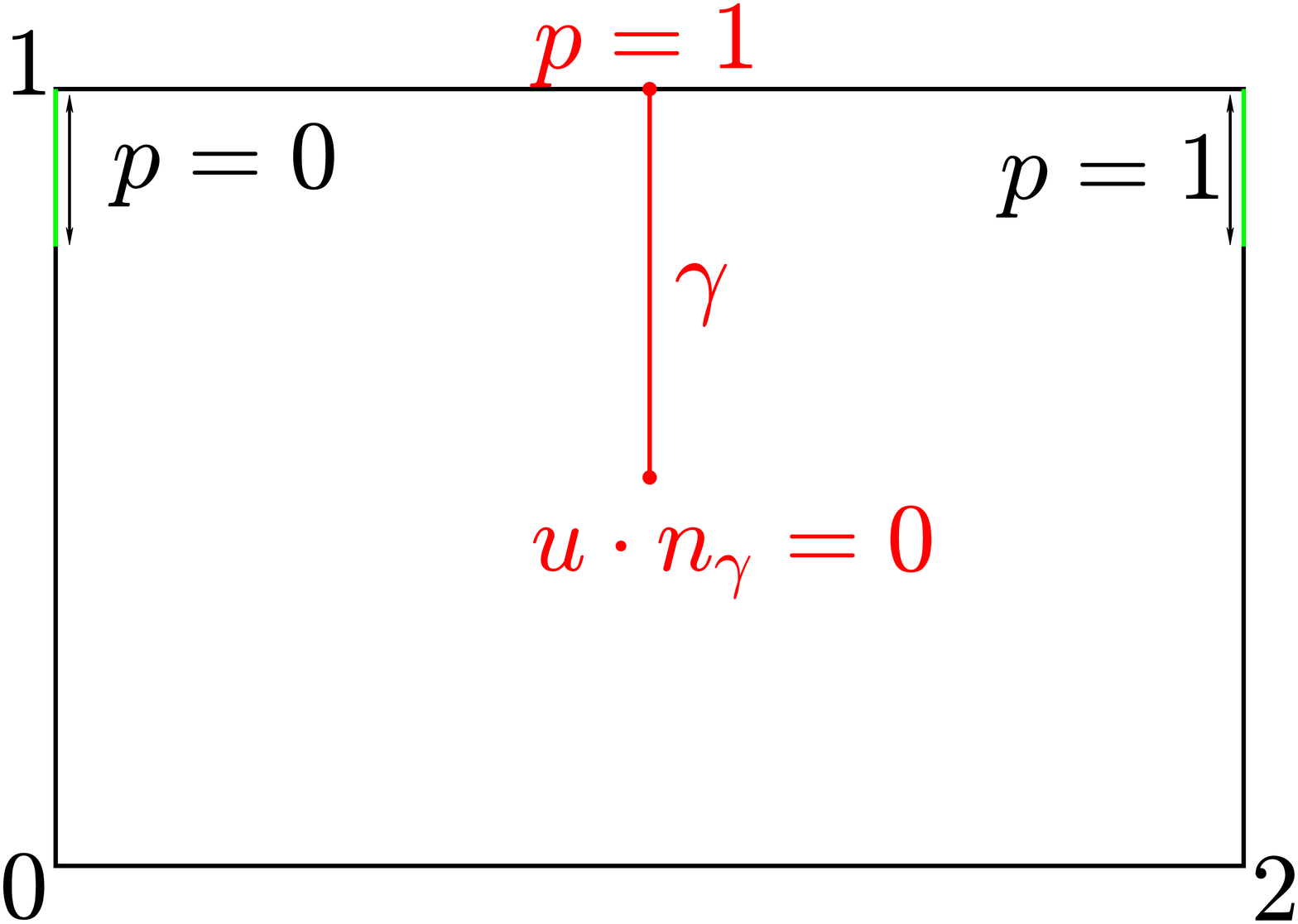}
  \label{fig:test1}
\end{minipage}%
\begin{minipage}{.5\textwidth}
  \centering
  \includegraphics[width=0.97\linewidth]{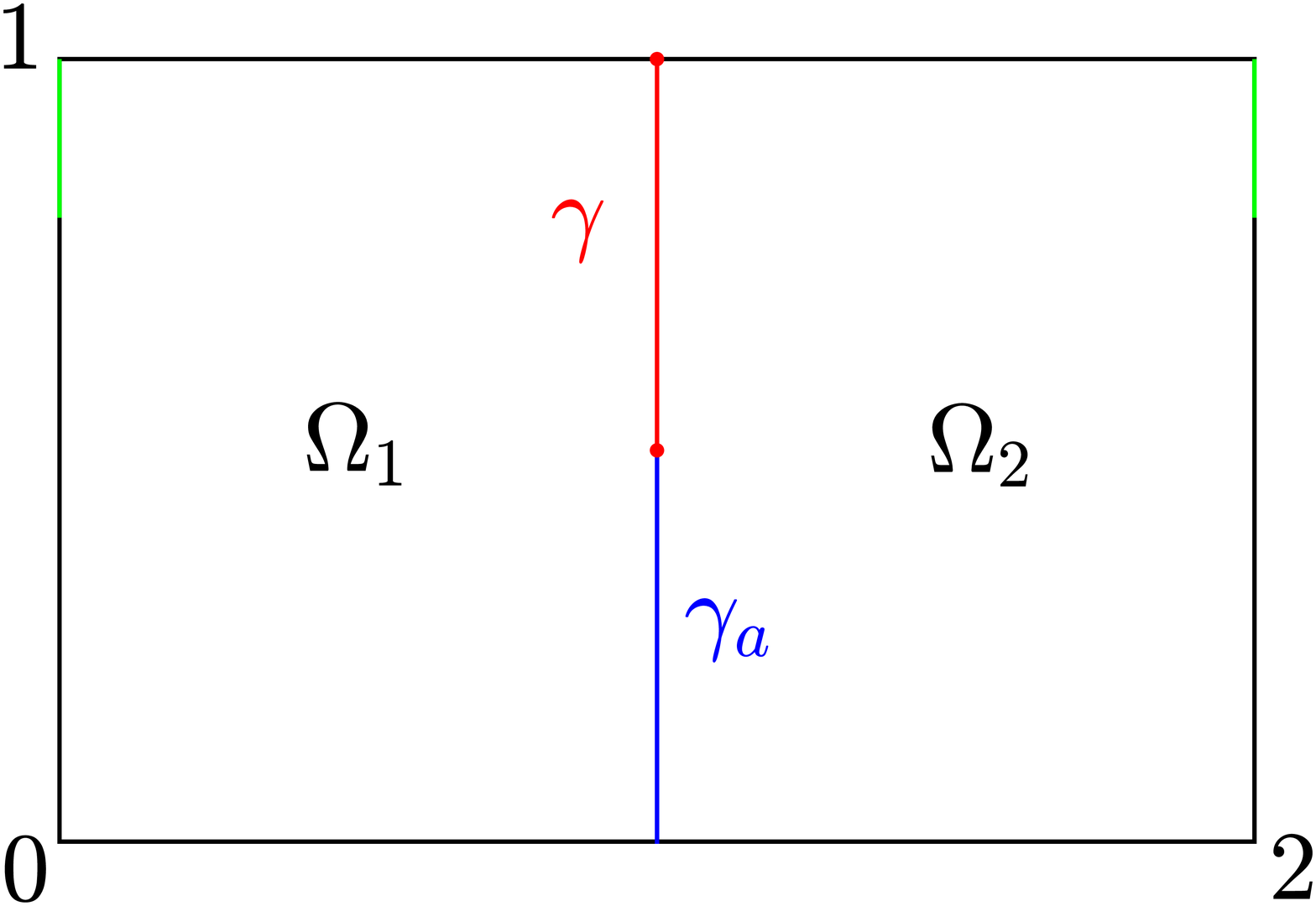}
  \label{fig:test2}
\end{minipage}
\caption{[Test case 2] (Left) Geometry and boundary conditions with immersed fracture $\gamma$. (Right) An artificial interface $\gamma_{a}$ is introduced to decompose the domain into two disjoint subdomains.}
\label{Immersed} \vspace{-0.4cm}
\end{figure}

We consider a test case adapted from \cite{4} where only one tip of the fracture is attached to the external boundary,  while the other tip is submerged inside the rock matrix as depicted in Figure~\ref{Immersed} (left).  A no-flow boundary condition is considered at the tip which is immersed inside the domain, while $p=1$ is imposed at the other tip.  Analysis of the steady-state flow problem with an immersed fracture can be found in \cite{4} and the references therein.  For the external boundary,  the pressure is prescribed on the upper fifth (length 0.2) of both lateral sides,  $p = 1$ on the right and $p = 0$ on the left, and a no flow condition is imposed on the rest of the boundary.  Note that we use the same physical parameters as in Test case 1.   

To apply global-in-time DD methods for this test case,  we first introduce an artificial interface $\gamma_a$ so that, together with the partially immersed fracture $\gamma$,  they form a single fracture $\Gamma$ separating the original domain into two disjoint subdomains (cf. Figure~\ref{Immersed} (right)).  Next,  suitable transmission conditions will be imposed on this new interface $\Gamma$.   On the fracture-interface $\gamma$,  we use the transmission conditions associated with the reduced fracture model (cf.  Equations~\eqref{transmission_1}-\eqref{transmission_2}).  Note that due to the presence of the immersed tip, we use a no-flow boundary condition at that tip,  instead of a Dirichlet condition as in Test case 1. 
On the artificial interface $\gamma_{a}$,  standard DD transmission conditions (representing the continuity of the pressure and normal flux) are imposed: \vspace{-0.2cm}
\begin{equation}
\label{normal_trans_conditions}
\begin{array}{c}
p_1 = p_2, \\
\textbf{\textit{u}}_1\cdot\textbf{\textit{n}}_1 + \textbf{\textit{u}}_2\cdot\textbf{\textit{n}}_2 =0,
\end{array} \qquad \text{on} \; \gamma_a \times (0, T). \vspace{-0.2cm}
\end{equation}
Global-in-time DD methods for parabolic equations in a domain without fractures have been well studied in \cite{24}, and thus will be omitted here.  The interface problem on $\Gamma \times (0,T)$ for each global-in-time DD method is then a combination of a problem on the interface-fracture (as derived in the previous sections for the non-immersed fracture case) and another one on the artificial interface as studied in \cite{24}.  For the latter,  we will also use preconditioners to enhance the convergence of the iterative algorithms. In particular, for GTP-Schur and GTF-Schur,  a time-dependent Neumann-Neumann preconditioner~\cite{24} is applied on the artificial interface,  while for GTD-Schur, a time-dependent Dirichlet-Dirichlet preconditioner is performed.  Combining the preconditioners on both the fracture-interface and artificial interface, we obtained the following methods: preconditioned GTP-Schur (with V-V preconditioner on the fracture-interface),  preconditioned GTD-Schur and preconditioned GTF-Schur.  These methods will be tested and compared with the performance of GTO-Schwarz.  Note that the transmission conditions for GTO-Schwarz on the artificial interface $\gamma \times (0,T)$ are Robin conditions with optimized parameters; more details can be found in \cite{24}. 

\begin{figure}[h!]
\centering
\vspace{-0.3cm}
\begin{minipage}{.5\textwidth}
  \centering
  \includegraphics[width=1\linewidth]{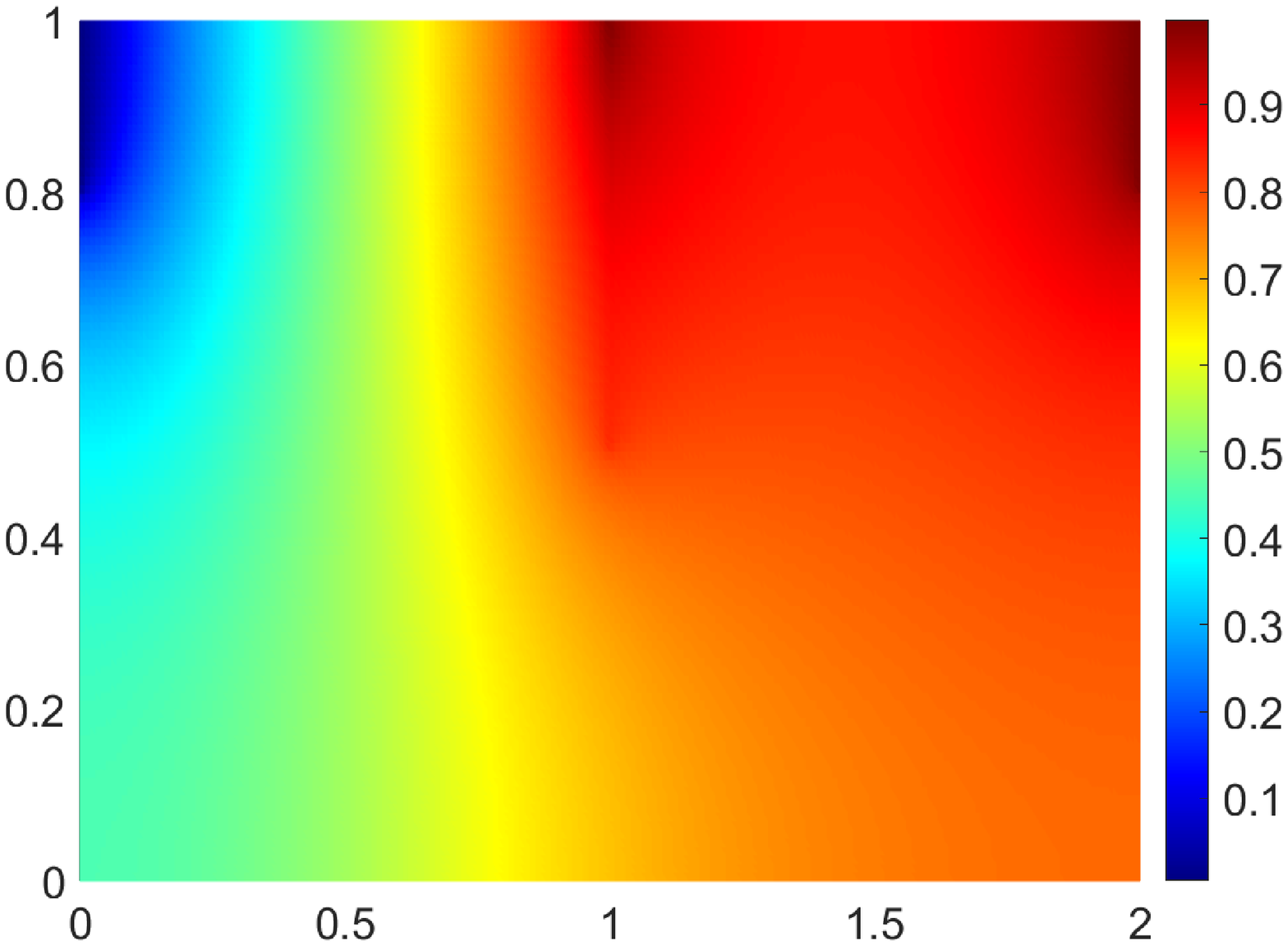}
  \label{fig:pressure_field_immersed}
\end{minipage}%
\begin{minipage}{.5\textwidth}
  \centering
  \includegraphics[width=1\linewidth]{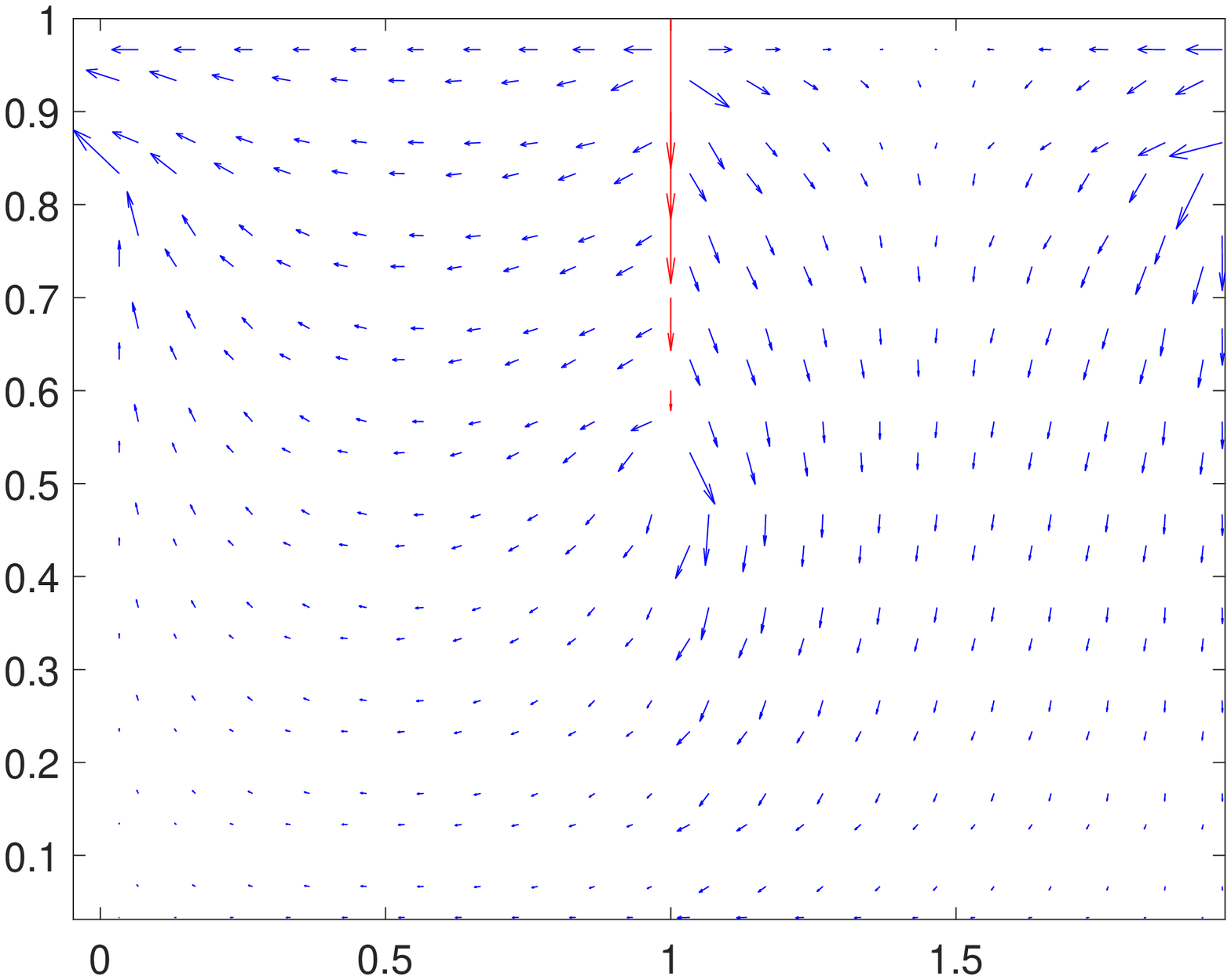}
  \label{fig:velocity_field_immersed}
\end{minipage}
\caption{[Test case 2] Pressure field (left) and velocity field (right) at the final time $T=1$.}
\label{Pressure_Vel_fields_Immersed} \vspace{-0.4cm}
\end{figure}

We first show the snapshots of pressure and velocity fields at the final time $T=1$ in Figure~\ref{Pressure_Vel_fields_Immersed}. The length of each arrow is proportional to the magnitude of the velocity and the red arrows represent the flow in the fracture. The length of the red arrows decreases as the flow travels toward the immersed tip since a no-flow boundary condition is imposed there.  As $\mathsf{k}_f \gg \mathsf{k}_i, \; i=1, 2$, the velocity in the fracture has larger magnitude than the one in the rock matrix.

Next, we present the numerical results for these methods when {\em conforming} time grids are used. We fix the spatial mesh $h=1/100$ and vary the time step sizes $\Delta{t}_i = \Delta{t}$ for $i=1, 2, \gamma$. The reference solution used in computing the errors are found on a fine time grid $\Delta{t}_{\text{ref}} = T/2000$ where $T=1$. Table~\ref{Err_Concentration_Velocity_Immersed} shows the relative $L^2$-errors computed from all methods for pressure and velocity. Similar to the non-immersed fracture case, we only have one table showing the errors each term since the approximate solutions obtained from all method are nearly the same. It can be observed that we still have first-order convergence in time for both pressure and velocity, even in the immersed fracture case which is more complicated. 

We report in Table~\ref{GMRES_iterations_Immersed_conforming} the number of subdomain solves needed to reach the errors obtained in Table~\ref{Err_Concentration_Velocity_Immersed}. It can be seen that the preconditioned GTD-Schur has the slowest convergence speed compared to the other methods, while the convergence speed of the preconditioned GTF-Schur method is acceptable. The preconditioned GTP-Schur method is still fast and comparable with the GTO-Schwarz method. Unlike Test case $1$, the preconditioned GTP-Schur method is the only method that has nearly the same speed as the GTO-Schwarz method.

\begin{table}[h]
\centering
\begin{minipage}{\textwidth}
\begin{tabular*}{\textwidth}{@{\extracolsep{\fill}}ccccccc@{\extracolsep{\fill}}}
\toprule%
& \multicolumn{3}{@{}c@{}}{Errors for pressure} & \multicolumn{3}{@{}c@{}}{Errors for velocity} \\\cmidrule{2-4}\cmidrule{5-7}%
$\Delta{t}$ & $\Omega_1$ & $\Omega_2$ & $\gamma$ & $\Omega_1$ & $\Omega_2$ & $\gamma$ \\
\midrule
T/4  & 7.85e-02 & 6.67e-02 & 2.55e-02 & 8.52e-02 & 2.88e-01 & 1.97e-01\\
&&&&& \\
T/8  & 4.15e-02 & 3.43e-02  & 1.27e-02  & 4.54e-02 & 1.49e-01 & 9.83e-02\\
& $[0.92]$ & $[0.96]$ & $[1.00]$ & $[0.91]$ & $[0.95]$ & $[1.00]$ \\
&&&&& \\
T/16 & 2.12e-02 & 1.73e-02 & 6.36e-03 & 2.34e-02 & 7.55e-02 & 4.93e-02  \\
& $[0.97]$ & $[0.99]$ & $[0.99]$ & $[0.96]$ & $[0.98]$ & $[0.99]$ \\
&&&&& \\
T/32 & 1.07e-02 & 8.65e-03 & 3.17e-03 & 1.18e-02 & 3.77e-02 & 2.46e-02 \\
& $[0.99]$ & $[1.00]$ & $[1.00]$ & $[0.99]$ & $[1.00]$ & $[1.00]$ \\
\botrule
\end{tabular*}
\caption{[Test case 2] Relative $L^2$ errors of the pressure and velocity with {\em conforming} time steps.  The corresponding convergence rates are shown in square brackets.}\label{Err_Concentration_Velocity_Immersed}
\end{minipage} \vspace{-0.3cm}
\end{table}
\begin{table}[h!]
\centering
\begin{minipage}{\textwidth}
\begin{tabular*}{\textwidth}{@{\extracolsep{\fill}}lccccc@{\extracolsep{\fill}}}
\toprule%
& $\Delta{t}$  & T/4 & T/8 & T/16 & T/32 \\%
Methods & &   &  &  & \\
\midrule
Preconditioned GTP-Schur && $16$&$16$& $16$& $18$ \\
Preconditioned GTD-Schur && $42$&$50$&$62$&$66$ \\
Preconditioned GTF-Schur &&$26$&$32$&$40$& $42$ \\
GTO-Schwarz &&$23$&$23$&$24$&$24$ \\
\botrule
\end{tabular*}
\caption{[Test case 2] Numbers of subdomain solves when {\em conforming} time steps are used; the tolerance for GMRES is set to be $10^{-8}$.}\label{GMRES_iterations_Immersed_conforming}
\end{minipage} \vspace{-0.3cm}
\end{table}

We next investigate the numerical performance of these methods with {\em nonconforming} time grids.  \Rv{For the preconditioned GTP-Schur and preconditioned GTD-Schur methods,  numerical results suggest that the initial guess  for GMRES needs to be rescaled to obtain accurate numerical solutions.  Such a rescaling is done in our numerical experiments by using the Heged\"{u}s formula (cf.  \cite[Chapter 5, Subsection 5.8.3]{Liesen2013}).}
The relative errors for pressure and velocity are presented in Table~\ref{Err_Pressure_Immersed_Nonconforming} and Table~\ref{Err_Velocity_Immersed_Nonconforming}. Similar to Test case 1, we impose the same large time step in the subdomains and a smaller one in the fracture: $\Delta{t}_1 = \Delta{t}_2 = \Delta{t} = 4\Delta{t}_{\gamma}$. We consider the same groups of errors as in Test case 1.  By comparing with Table~\ref{Err_Concentration_Velocity_Immersed}, we can see that the fine time grids in the fracture do not affect the errors in the fracture for both pressure and velocity observed from Group 1, that is, we still obtain the same errors as when we only have coarse time grids in the subdomains and in the fracture.  On the contrary, such errors provided by Group $2$ are smaller, and closer to the ones obtained when we apply the same fine time grids in the subdomains and the fracture. These behaviors are as expected as explained in Test case 1.

\begin{table}[h!]
\centering
\begin{minipage}{\textwidth}
\begin{tabular*}{\textwidth}{@{\extracolsep{\fill}}cccccccc@{\extracolsep{\fill}}}
\toprule%
& & \multicolumn{3}{@{}c@{}}{{Preconditioned GTP-Schur}} & \multicolumn{3}{@{}c@{}}{Preconditioned GTD-Schur}  \\ 
& & \multicolumn{3}{@{}c@{}}{{GTO-Schwarz}} & \multicolumn{3}{@{}c@{}}{Preconditioned GTF-Schur}  \\ 
 \cmidrule{3-5} \cmidrule{6-8}
$\Delta{t}$ & $\Delta{t}_{\gamma}$ & $\Omega_1$ & $\Omega_2$ & $\gamma$ & $\Omega_1$ & $\Omega_2$ & $\gamma$\\
\midrule
T/4  & T/16 & 7.85e-02 & 6.67e-02 & 2.55e-02 & 7.61e-02 & 6.51e-02 & 1.50e-02 \\
&&& \\
T/8 & T/32 & 4.15e-02 & 3.43e-02  & 1.27e-02 & 3.98e-02 & 3.33e-02 & 7.32e-03 \\
& & $[0.92]$ & $[0.96]$ & $[1.00]$ & $[0.93]$ & $[0.97]$ & $[1.03]$  \\
&&& \\
T/16 & T/64 & 2.12e-02 & 1.73e-02&  6.36e-03 & 2.04e-02 & 1.67e-02 & 3.58e-03\\
& & $[0.97]$ & $[0.99]$ & $[0.99]$ &  $[0.96]$ & $[0.99]$ & $[1.03]$  \\
&&& \\
T/32 & T/128 & 1.07e-02 & 8.65e-03&  3.17e-03 & 1.02e-02 & 8.33e-03 & 1.75e-03\\
& & $[0.99]$ & $[1.00]$ & $[1.00]$ & $[1.00]$ & $[1.00]$ & $[1.03]$  \\
\botrule
\end{tabular*}
\caption{[Test case 2] Relative $L^2$ errors of the {\em pressure} with nonconforming time grids. The corresponding convergence rates are shown in square brackets. }\label{Err_Pressure_Immersed_Nonconforming}
\end{minipage} \vspace{-0.4cm}
\end{table}
\begin{table}[h!]
\centering
\begin{minipage}{\textwidth}
\begin{tabular*}{\textwidth}{@{\extracolsep{\fill}}cccccccc@{\extracolsep{\fill}}}
\toprule%
& & \multicolumn{3}{@{}c@{}}{{Preconditioned GTP-Schur}} & \multicolumn{3}{@{}c@{}}{Preconditioned GTD-Schur}  \\ 
& & \multicolumn{3}{@{}c@{}}{{GTO-Schwarz}} & \multicolumn{3}{@{}c@{}}{Preconditioned GTF-Schur}  \\ 
 \cmidrule{3-5} \cmidrule{6-8}
$\Delta{t}$ & $\Delta{t}_{\gamma}$ & $\Omega_1$ & $\Omega_2$ & $\gamma$ & $\Omega_1$ & $\Omega_2$ & $\gamma$\\
\midrule
T/4  & T/16 & 8.51e-02 & 2.88e-01 & 1.97e-01 & 8.38e-02 & 2.86e-01 & 1.18e-01 \\
&&& \\
T/8  & T/32 & 4.54e-02 & 1.49e-01  & 9.83e-02 & 4.41e-02 & 1.47e-01 & 5.73e-02 \\
& & $[0.91]$ & $[0.95]$ & $[1.00]$ & $[0.93]$ & $[0.96]$ & $[1.04]$ \\
&&& \\
T/16 & T/64 & 2.34e-02 & 7.55e-02& 4.93e-02 & 2.26e-02 & 7.36e-02 & 2.80e-02\\
& & $[0.96]$ & $[0.98]$ & $[0.99]$ & $[0.96]$ & $[0.99]$ & $[1.03]$ \\
&&& \\
T/32 & T/128 & 1.18e-02 & 3.77e-02&  2.46-02 & 1.13e-02 & 3.66e-02 & 1.37e-02\\
& & $[0.99]$ & $[1.00]$ & $[1.00]$ & $[1.00]$ & $[1.01]$ & $[1.03]$\\
\botrule
\end{tabular*}
\caption{[Test case 2] Relative $L^2$ errors of the {\em velocity} with nonconforming time grids. The corresponding convergence rates are shown in square brackets.}\label{Err_Velocity_Immersed_Nonconforming} 
\end{minipage} \vspace{-0.4cm}
\end{table}

Finally, we present the number of subdomain solves for each method to reach the relative residual smaller than $10^{-8}$ to analyze their convergent behaviors. These numbers are shown in Table~\ref{GMRES_iterations_Immersed_nonconforming}. It can be seen that we obtain nearly the same numbers as those in Table~\ref{GMRES_iterations_Immersed_conforming}. Hence, as in Test case 1, these methods are applicable under nonconforming time discretizations. From what we have observed so far, Test case 2 is more challenging than Test case 1, which can be seen in the increasing of the subdomain solves. However, the preconditioned GTF-Schur still shows its efficiency as it has relatively fast convergence speed and preserves the accuracy in time when we have different time steps in the fracture and in the subdomains.

\begin{table}[h!]
\centering
\setlength{\extrarowheight}{4pt}
\begin{tabular*}{\textwidth}{@{\extracolsep{\fill}}lccccc@{\extracolsep{\fill}}}
\toprule%
& $\Delta{t}_1 = \Delta{t}_2$ & T/4 & T/8 & T/16 & T/32 \\
& $\Delta{t}_{\gamma}$ & T/16 & T/32 & T/64 & T/128 \\
\midrule
Preconditioned GTP-Schur  && $16$&$14$& $14$& $14$ \\
Preconditioned GTD-Schur &&$42$&$50$&$60$&$66$ \\
Preconditioned GTF-Schur &&$26$&$32$&$40$& $44$ \\
GTO-Schwarz &&$23$&$24$&$24$&$24$ \\
\botrule
\end{tabular*}
\caption{[Test case 2] Numbers of subdomain solves when {\em nonconforming} time steps are used; the tolerance for GMRES is set to be $10^{-8}$. } 
\label{GMRES_iterations_Immersed_nonconforming} 
\end{table}
\section*{Conclusion}
In this work, three global-in-time DD methods, namely GTP-Schur, GTD-Schur and GTF-Schur, have been developed for a reduced fracture model of compressible flow problems,  in which different time steps can be used in the fracture and in the matrix.  Efficient preconditioners have been derived for GTP-Schur and GTD-Schur to enhance the convergence of the iterative algorithms. Importantly,  a new method, GTF-Schur, is proposed; this method is typical to the reduced fracture model and requires no preconditioner. 
Numerical experiments with different types of fractures have been carried out to investigate the performance of the proposed methods on conforming and nonconforming time grids.  The obtained results suggest that GTF-Schur is the most efficient method as it converges fast without preconditioning while preserving the accuracy in time in the fracture when smaller time steps are used in the fracture and larger ones in the rock matrix. 
Our ongoing work is to extend these methods to solve the advection-diffusion problem with operator splitting,  in which the advection is treated explicitly and the diffusion implicitly.  Such an approach gives satisfactory results when advection is mild.  For strongly advection-dominated problems,  we will use mixed-hybrid finite element method proposed in \cite{43, 8}, and develop corresponding global-in-time DD methods based on both physical and optimized transmission conditions.  

\bmhead{Acknowledgments} This work is partially supported by the US National Science Foundation under grant numbers DMS-1912626 and DMS-2041884.

\bibliography{sn-bibliography}


\end{document}